\documentclass[a4paper]{amsart}
\usepackage{amsmath}
\usepackage{amsfonts}
\usepackage{amssymb}
\usepackage{latexsym}
\usepackage{epsfig}

\def\co{\colon\thinspace}
\newcommand{\begriff}[1]{\textbf{#1}}

\newcommand{\p}{\mathbb S}
\renewcommand{\d}{\partial}

\newcommand{\C}{\mathcal C}

\newtheorem{theorem}{Theorem}

\newtheorem{lemma}{Lemma}
\newtheorem{corollary}{Corollary}

\theoremstyle{definition}
\newtheorem{definition}{Definition}

\newtheorem{example}{Example}

\theoremstyle{remark}

\begin{document}
\bibliographystyle{alpha}

\title{{\bf 
A TOPOLOGICAL REPRESENTATION THEOREM\\
 FOR ORIENTED MATROIDS\\}
}
\setcounter{footnote}{-1}
\author{
J\"urgen Bokowski \hspace{4mm} Simon King \hspace{4mm}  Susanne Mock
%\thanks{
%Department of Mathematics,
%        Darmstadt University of Technology, Schlossgartenstr.~7, \\
%        D-64289 Darmstadt, Germany.
%bokowski@mathematik.tu-darmstadt.de, \\
%king@mathematik.tu-darmstadt.de,\\
%mock@mathematik.tu-darmstadt.de
%}
\hspace{4mm}
Ileana Streinu
%\thanks{Dept. of Computer Science, Smith College, Northampton, MA 01063, USA,
%streinu@cs.smith.edu}
} %end author

\date{\today}

\begin{abstract}
We present a new direct proof of 
a topological representation theorem for oriented matroids
in the general rank case. Our proof is based on an earlier rank 3 version.
It uses hyperline sequences and the generalized Sch{\"o}nflies theorem.
%new
As an application, we show that one can read off oriented matroids
from arrangements of embedded spheres of codimension one, even if wild spheres are involved.
\end{abstract}

\maketitle


\section{Introduction}
\label{introduction}

When studying vector configurations or central hyperplane configurations,
point sets on a sphere or great hypersphere arrangements, vector spaces
or their duals, points on grassmannians, polytopes and their corresponding
cellular decompositions in projective space, etc., an abstraction 
of an important equivalence class of matrices often plays a central
role: an oriented matroid.
The theory of oriented matroids (see \cite{bjorner}) 
provides us with a multitude of definitions
for an oriented matroid that can be viewed as
reflecting the variety of objects that a matrix can represent.

These definitions via 
circuit or cocircuit axioms (\cite{bjorner}, p.~103), sphere systems
(\cite {bjorner}, p.~227), Gra\ss mann Pl\"ucker relations (chirotope
axioms) (\cite {bjorner}, p.~126, p.~138, \cite{GN}), hull systems
(\cite{knuth}), to mention just a few of them, differ a lot with
respect to their motivational aspects, their algorithmical efficiency
or their relation to the actual application.  Each definition provides
in general an additional insight for the motivating problem. In the
research monograph on oriented matroids, \cite{bjorner}, three
chapters are devoted to axiomatics concerning oriented matroids and to
the topological representation theorem for oriented matroids.

This central theorem in the theory of oriented matroids
due to Lawrence shows the equivalence of 
oriented matroids defined via sphere system axioms with
oriented matroids defined via, say, the circuit axioms.
This remarkable result asserts that each oriented matroid 
has a topological representation as an oriented pseudosphere arrangement,
even a piecewise-linear one, cf. Edmonds and Mandel \cite{edmonds}.
Other authors 
(\cite{bjorner}, \cite{hochstattler}) have later simplified or complemented 
the original proof, but all use fundamentally the same approach: the
face lattice (tope) formalism for oriented matroids and a shelling order to
carry through the construction.

Finding a reasonably direct proof in rank $3$, 
one that would rely on the structural
simplicity of the planar case,
has been posed as an open problem in the research monograph \cite{bjorner}
(Exercise 6.3). 
In \cite{bokowskimockstreinu} such a proof in the rank 3 case was given.
Unlike the previous ones, this was based on 
\emph{hyperline sequences}, an equivalent axiomatization for oriented 
matroids which is particularly natural in rank $3$. 
In this article we generalize this proof to the arbitrary rank case.
The proof is inductive, direct and uses only one advanced result from
topology, the generalized Sch\"onflies theorem.

The motivation for introducing hyperline sequences in rank 3 was
presented in detail in our previous paper~\cite{bokowskimockstreinu}.
Here we use the generalization to arbitrary rank, formally defined in
Section~\ref{sec:HLS}.  The motivation is similar. Here is a brief
account of it. 
We use an index set $E_n = \{1,\dots,n\}$ with its standard total
order.  We consider a labeled finite point set $X:=\{x_1,x_2,\ldots,
x_n\},$ $x_i \in \mathbb R^{r-1}$, ordered according to their labels.
We assume that the point set $X$ affinely spans $\mathbb R^{r-1}$. 
We endow the affine Euclidean space $\mathbb R^{r-1}$ with an orientation, 
i.e., for any $r$--element subset of affinely independent points, 
we can tell whether the ordering of their indices 
defines a positive, or a negative, oriented simplex, respectively.
Every $(r-1)$--element subset of affine independent points of $X$ 
defines via the ordering of their indices an oriented hyperplane, 
i.e., an affine hull of codimension $1$ together with an orientation.
Furthermore, every $(r-2)$--element subset of affine independent points of $X$ 
defines via the ordering of their indices an oriented hyperline,
i.e., an affine hull of codimension 2 together with an orientation of it. 
Although this affine model can be described further and 
the concept we are going to describe
is useful for affine applications as well, 
in what follows we prefer to use the notational advantage that occurs 
when embedding this concept in the usual way in projective space. 
First consider a bijective map between our former 
oriented $(r-1)$--dimensional Euclidean space $\mathbb R^{r-1}$ 
and a non-central oriented hyperplane of the oriented space $\mathbb R^{r}$. 
The former points $x_i$ then become non-zero vectors $v_i$ in $\mathbb R^{r}$,
$i \in E_n$, the former oriented hyperplanes become
oriented central hyperplanes in $\mathbb R^{r}$, 
and the former oriented hyperlines, defined via an index set 
$B \subset E_n$ with $|B|=r-2$, become 
oriented central hyperlines (subspaces of codimension 2 in $\mathbb R^{r}$ 
together with an orientation). For what follows we can also assume 
that we have started with an ordered set of vectors 
$V=\{v_1,v_2,\ldots, v_n\}$ in our oriented space $\mathbb R^{r}$ 
that has not necessarily an affine pre-image. 

The oriented central hyperline with index set $B$ is an $\mathbb
R^{r-2}$, and the restriction of  $V$ to all elements in that hyperline yields an
arrangement $Y$ of non-zero vectors in $\mathbb R^{r-2}$. 
We rotate an additional central oriented hyperplane in $\mathbb R^{r}$ 
around one of our oriented central hyperlines with index set $B$ 
in the well defined positive sense, i.e., the ordered set of vectors with
index set $B$  together with two adjacent consecutive outer normal
vectors of the hyperplane form a positive basis of our oriented space 
$\mathbb R^{r}$. 
An orthogonal projection $\pi_B$ along the oriented hyperline 
onto its corresponding two-dimensional orthogonal oriented space 
$\pi_B(\mathbb R^r)$  
(that forms together with the oriented hyperline a positive orientation 
of the whole space) shows the image of the additional central 
oriented hyperplane as a rotating (mathematically positive) 
central oriented line in the plane.
Precisely all vectors $v_i$ 
that are not in the hyperline defined by $B$ 
appear after the projection as non-zero vectors $\pi_B(v_i).$
Their indices form a periodic circular sequence around the origin
when we pursue the incidences of the vectors with the rotating line.
We store the cyclic order of incidences in a cyclic order of sets,
where we write an index $i$ when the vector $v_i$ is consistent with the
orientation of the oriented line, and $\overline{\imath}$
otherwise. We refer to this notation as \begriff{signed indices}.

When an incidence position, say position $m$, of the rotating 
oriented line with a vector occurs, we write the signed indices of all vectors 
that are incident with the oriented line in a set $Z^m$. Note that the
period of the cyclic sequence is even, since there
are two positions of the rotating oriented line in which it is  incident
with the vector $\pi_B(v_i)$, corresponding to the signed indices $i$ and
$\overline{\imath}$.
We obtain a well defined oriented circular sequence $Z^0,\dots,Z^{2k-1}$ 
(\begriff{oriented cycle}) of sets of signed indices, where
after half the period there is a sign reversal for each element.
The combinatorics of the hyperline is encoded in the arrangement $Y$ 
with index set $B$ together with the oriented circular sequence of incidences 
with the rotating line,
$(Y|Z^0,\ldots,Z^{2k-1})$,
to which we also refer as a \begriff{hyperline}.
We obtain a \begriff {hyperline sequence} when we write
down all hyperlines arising in the vector arrangement.

This model can be generalized in two ways.
The combinatorial abstraction leads on the one hand 
to oriented matroids characterized
as hyperline sequences and on the other hand it leads 
to topological equivalence classes of arrangements of pseudospheres.
%new
We give a new proof of a one-to-one correspondence of oriented matroids
and classes of arrangements of pseudospheres.
We  further generalize arrangements of pseudospheres,
and we show that one can still read off a hyperline sequence.

The paper is organised as follows.  We define hyperline sequences as a
combinatorial abstraction of vector configurations in 
Section~\ref{sec:HLS}. 
We recall the chirotope concept in Section~\ref{sec:chirotope} and we
show the equivalence of these two concepts in
Theorem~\ref{thm:HLS-Chi}.  Furthermore, we discuss the concepts of
deletion and contraction for both settings.
We introduce the topological representation of oriented matroids via
arrangements of oriented pseudospheres in Section~\ref{sec:PSA} with
the corresponding concepts of deletion and contraction. In the proof
of Theorem~\ref{thm:cells} we use the generalized Sch\"onflies theorem
to show the cell structure of an arrangement of oriented
pseudospheres.  In Theorem~\ref{thm:PSA} we replace our two axioms for
an arrangement of oriented pseudospheres with a single one.
Section~\ref{sec:X(A)} deals with the easier part of the main
representation theorem: we obtain chirotopes and hyperline sequences
from the topological representation.
We complete the proof of the one-to-one correspondence of hyperline
sequences with their topological representation by induction.
Section~\ref{sec:FLbase} is devoted to the base case and
Section~\ref{sec:FLgeneral} contains the essential part.
We finally discuss in Section~\ref{sec:wild} the wild arrangement case
as an easy supplement of our approach.

The whole proof is based heavily on the ideas that have been worked
out already in \cite{bokowskimockstreinu} by the first and the last
two authors in the rank~$3$ case.  But the experience with respect to
topological arguments of the second author was decisive to arrive at
our final version of the proof in the general rank case.  The last
section and many improvements of the proof compared with the rank~3
version are due to him. For instance, the uniform and non-uniform
cases form no longer separate parts within the proof.

%************************************************************************

\section{Hyperline sequences}
\label{sec:HLS}

Our aim in this section is to introduce the notion of hyperline
sequences.
The geometric motivation of our definition comes from vector arrangements in
Euclidian space as explained in the introduction.

Let $(E,<)$ be a finite totally ordered set. 
Let $\overline E = \{\overline e | e\in E\}$ be a copy of $E$. 
The set $\mathbf E$ of \begriff{signed indices} is defined as the
disjoint union of $E$ and $\overline E$. 
By extending the map $e\mapsto \overline e$ to $\overline e\mapsto
\overline{\overline e} = e$ for $e\in E$, we get an involution
on $\mathbf E$. 
We define $e^*= \overline{e}^* = e$. 
For $X\subset \mathbf E$, define $\overline X= \{\overline x|\;
x\in X\}$ and $X^* = \{x^* |\; x\in X\}$. 

An \begriff{oriented $d$--simplex} in $E$ is a 
$(d+1)$--tuple $\sigma = [x_1,\dots,x_{d+1}]$ of elements of $\mathbf E$,
so that $x_1^*,\dots, x_{d+1}^*$ are pairwise distinct. 
Let an equivalence relation $\sim$ on oriented $d$--simplices in $E$ be
generated by 
$[x_1,\dots,x_{d+1}]\sim [x_1,\dots, x_{i-1},\overline{x_{i+1}},x_i,
x_{i+2},\dots, x_{d+1}]$, for $i=1,\dots, d$. 
As usual, any oriented $d$--simplex is equivalent to one of
the form $[ e_1,\dots,  e_{d+1}]$ or $[ e_1,\dots,  e_{d},\overline{
  e_{d+1}}]$, with elements $e_1<e_2<\dots < e_{d+1}$ of $E$. 
Define $-[x_1,\dots,x_{d+1}] = [x_1,\dots,x_d,\overline{x_{d+1}}]$.
If $\phi\co \mathbf E\to \mathbb R^d$ is a
map with $\phi(\overline e) = -\phi(e)$ for all $e\in E,$ then
$\phi(x_1),\dots, \break \phi(x_{d+1})$ are the vertices of a simplex in
$\mathbb R^d$. Note that this simplex might be degenerate.

In the following iterative definition of hyperline sequences we denote
with $$C_m=\big(\{0,1,\dots, {m-1}\},+\big)$$
the cyclic group of
order $m$.
\begin{definition}[Rank $1$]
  A \begriff{hyperline sequence} $X$ over $E(X)\subset E$ of rank $1$
  is a non-empty subset $X\subset E(X)\cup \overline{E(X)}$ so that
  $|X| = |X^*|$ and $X^*=E(X)$. 
\end{definition}

The oriented simplex $[x]$ is by definition a
\begriff{positively oriented base} of $X$ for $x\in X$. 
We define $-X = \overline X$.

\begin{definition}[Rank $2$]
  Let $k\in \mathbb N$, $k\ge 2$.  A \begriff{hyperline sequence} $X$
  of rank $2$ over $E(X)\subset E$ is a map from $C_{2k}$ to hyperline
  sequences of rank one, $a\mapsto X^a$, so that $X^{a+k} =
  -X^a$ for all $a\in C_{2k}$, and $E(X)\cup\overline{E(X)}$
  is a disjoint union of $X^0,\dots, X^{2k-1}$.  
\end{definition}

We refer to $X^0,\dots, X^{2k-1}$ as the \begriff{atoms} of $X$ and
to $2k$ as the \begriff{period length} of $X$.
We say that $e\in E(X)$ is \begriff{incident} to an atom $X^a$ of $X$
if $e\in (X^a)^*$. Let $x_1,x_2\in E(X)\cup \overline{E(X)}$ so that
$x_1^*$ and $x_2^*$ are not incident to a single atom of $X$, and $X$
induces the cyclic order
$(x_1,x_2,\overline{x_1},\overline{x_2})$. Then, the oriented simplex
$[x_1,x_2]$ is by definition a \begriff{positively oriented base} of
$X$.  
We define the hyperline sequence $-X$ over $E(-X)=E(X)$ of rank $2$ as
the map $a\mapsto (-X)^a = X^{-a}$ for $a\in C_{2k}$.  

A hyperline sequence $X$ of rank $2$ is determined by the sequence
$(X^0,\dots, X^{2k-1})$ of atoms.  We define that two hyperline
sequences $X_1$ and $X_2$ of rank $2$ are equal, $X_1=X_2$, if
$E(X_1)=E(X_2)$, the number $2k$ of atoms coincides, and $X_1$ is
obtained from $X_2$ by a shift, i.e., there is an $s\in C_{2k}$ with
$X_1^{a+s} = X_2^a$ for all $a\in C_{2k}$.

We prepare the axioms for hyperline 
sequences of rank $r>2$ with the following definitions. 
Let $X$ be a set of pairs $(Y|Z)$, where $Y$ is a hyperline sequence 
of rank $r-2$ and $Z$ is a hyperline sequence of rank $2$.
A \begriff{positively oriented base} of $X$ in $(Y|Z)\in X$ is an
oriented simplex $[x_1,\dots,x_{r}]$ in $E(X)$, where $[x_1,\dots
x_{r-2}]$ is a positively oriented base of $Y$ and $[x_{r-1},x_r]$
is a positively oriented base of $Z$.
We define $-X = \{ (Y|-Z) \;|\; (Y|Z)\in X\}$. The elements of
  $X$ are called \begriff{hyperlines}. An \begriff{atom} of $X$ in a
  hyperline $(Y|Z)\in X$ is the pair $(Y|Z^a)$, where $Z^a$ is an atom
  of $Z$.

\begin{definition}[Rank $r>2$]
  Let $X \not=\emptyset$ be a set whose elements are pairs $(Y|Z)$,
  where $Y$ is a hyperline sequence of rank $r-2$ and $Z$ is a
  hyperline sequence of rank $2$.
  The set $X$ is a \begriff{hyperline sequence} of rank $r>2$ over
  $E(X)\subset E$ if it satisfies the following axioms.
  \begin{itemize}
  \item[(H1)] $E(X)$ is a disjoint union of $E(Y)$ and $E(Z)$, for all
    $(Y|Z)\in X$.
  \item[(H2)] Let $(Y_1|Z_1), (Y_2|Z_2)\in X$ and let $[x_1,\dots
    x_{r-2}]$ be a positively oriented base of $Y_1$.
    If $\{x_1^*,\dots, x_{r-2}^*\}\subset E(Y_2)$ then 
    $(Y_1|Z_1) = (Y_2|Z_2)$ or $(Y_1|Z_1) = (-Y_2|-Z_2)$.
  \item[(H3)] For all positively oriented bases $[x_1,\dots, x_r]$ and
    $[y_1,\dots,y_r]$ of $X$, there is some $j\in\{1,\dots,r\}$ so
    that $[x_1,\dots, x_{r-1},y_j]$ or $[x_1,\dots,
    x_{r-1},\overline{y_j}]$ is a positively oriented base of~$X$.
  \item[(H4)] For any positively oriented base $[x_1,\dots, x_r]$ of
    $X$, $[x_1,\dots,x_{r-3},\overline{x_{r-1}},x_{r-2},x_r]$ is a
    positively oriented base of $X$.
  \end{itemize}
\end{definition}

We call an oriented simplex $\sigma$ a negatively oriented base of a
hyperline sequence $X$ if $-\sigma$ is a positively oriented base of
$X$.

We connect these axioms to the geometric motivation exposed in the
introduction. Let $V=\{v_1,\dots, v_n\}\subset \mathbb R^{r-1}$ be
non-zero vectors spanning $\mathbb R^{r}$. Let $H\subset \mathbb
R^{r}$ be a hyperline, i.e. a subspace of dimension $r-2$ spanned by
elements of $V$. Let $H^\perp$ be the orthogonal complement of $H$ in
$\mathbb R^{r}$.  Then $V\cap H$ is a set of non-zero vectors in
$\mathbb R^{r-2}$, corresponding to the hyperline sequence $Y$ in a
hyperline $(Y|Z)$.
The image of $V\setminus H$ under the orthogonal projection onto
$H^\perp$ yields a set of non-zero vectors in $\mathbb R^2$
corresponding to the term $Z$ in $(Y|Z)$.
Axiom~(H1) simply means that $V$ is a disjoint union of $V\cap H$ and
$V\setminus H$.  Axiom~(H2) corresponds to the fact that $H$ is
uniquely determined by any $(r-2)$--tuple of elements of $V\cap H$ in
general position.
Axiom~(H3) is the Steinitz--McLane exchange lemma, stating that one
can replace any vector in a base by some vector of any other base.
Axiom~(H4) ensures that if $[x_1,\dots ,x_r]$ is a positively
oriented base of $X$ than so is any oriented simplex that is equivalent
to $[x_1,\dots, x_r]$; this is a part of Theorem~\ref{thm:HLS-Chi} below.
Axiom~(H4) is related to the 
``consistent abstract sign of
determinant'' in~\cite{bokowskimockstreinu}. It means that if $r$
points span an $(r-1)$--simplex, then any subset of $r-2$ points spans
a hyperline, and the orientation of the $(r-1)$--simplex does not
depend on the hyperline on which we consider the $r$ points.

%******************************************************************

\section{Chirotopes}
\label{sec:chirotope}

We recall in this section the chirotope axioms for oriented matroids
(see \cite {bjorner}, p.~126,p.~138, \cite{GN}). 
Let $(E,<)$ be as in the preceding section.  We
denote by $\Delta_d(E)$ the set of all oriented $d$--simplices in $E$.

\begin{definition}
  A \begriff{chirotope} $\chi$ of rank $r$ over $E$ is a map
  $\Delta_{r-1}(E) \to \{-1,0,+1\}$, so that the following holds.
  \begin{itemize}
  \item[(C1)] For any $e_1\in E$, there are $e_2,e_3,\dots,e_{r}\in E$
    with $\chi([ e_1,\dots, e_r])\not=0$.
  \item[(C2)] For any $\sigma\in \Delta_{r-1}(E)$ holds $\chi(-\sigma)
    = -\chi(\sigma)$, and if $\sigma\sim \tau$ then $\chi(\sigma) =
    \chi(\tau)$.
  \item[(C3)] If $\chi([x_1,\dots,x_r])\not=0$ and
    $\chi([y_1,\dots,y_r])\not=0$ then there is some
    $i\in\{1,\dots,r\}$ with $\chi([x_1,\dots,x_{r-1}, y_i])\not=0$.
  \item[(C4)] If
    $x_1,x_2,\dots, x_r, y_1,y_2\in \mathbf E$ so that
    \begin{eqnarray*}
        \chi([x_1,\dots,x_{r-2},y_1,x_r])\cdot \chi([x_1,\dots,x_{r-1},y_2])
        &\ge& 0\text{ and}\\
        \chi([x_1,\dots,x_{r-2},y_2,x_r]) \cdot 
        \chi([x_1,\dots,x_{r-1},\overline{y_1}])&\ge& 0,
    \end{eqnarray*}
 then $\chi([x_1,\dots,x_r])\cdot\chi([x_1,\dots,x_{r-2},y_1,y_2])\ge 0$.
  \end{itemize}
\end{definition}

\noindent
A simplex $\sigma\in\Delta_{r-1}(E)$ is a \begriff{positively oriented
  base} of $\chi$ if $\chi(\sigma) = +1$.  
By Axiom~(C2), $\chi$ is completely described by the set of its positively 
oriented bases.

In order to give Axioms~(C1)--(C4) a
geometric meaning, we show how to construct a chirotope of rank $r$ from a 
vector configuration in $\mathbb R^r$.
Let $V=\{v_1,\dots v_n\}\subset \mathbb R^{r}$ be a set of non-zero vectors
spanning $\mathbb R^{r}$. Let $\phi\co \mathbf E\to \mathbb R^{r}$ be
defined by $\phi(e) = v_e$ and $\phi(\overline e) = -v_e$ for all $e\in
E$.  We get a map $\chi\co \Delta_{r-1}(E) \to \{-1,0,+1\}$ by setting
\begin{itemize}
\item $\chi([x_1,\dots,x_r])=+1$ if the determinant
  $[\phi(x_1),\dots,\phi(x_r)]$ is positive 

\noindent  (i.e. if
  $\phi(x_1),\dots,\phi(x_r)$ is a positive base of $\mathbb R^{r}$),
\item $\chi([x_1,\dots,x_r]) = -1$ if $[\phi(x_1),\dots,\phi(x_r)]$
  is negative, and
\item $\chi([x_1,\dots,x_r]) =0$ if $\phi(x_1),\dots,\phi(x_r)$ are
  linearly dependent.
\end{itemize}

We show that $\chi$ is a chirotope. Axiom~(C1) holds since $V$ spans
$\mathbb R^{r}$. Axiom~(C2) follows from the symmetry of
determinants. 
Axiom~(C3) is the Steinitz--McLane exchange lemma, as Axiom~(H3)
above. Axiom~(C4) sais that the signs of oriented bases as defined by
$\chi$ do not contradict the three summand Gra{\ss}mann--Pl\"ucker
relation.  

Note that we do distinguish between $\chi$ and $-\chi$,
although this is not usual in the literature. Our reason comes from
geometry. We believe if one deals with oriented objects (vectors,
oriented hyperspheres, etc.) in $\mathbb R^d$ or $\p^d$, then it is
consequent to distinguish not only the two different orientations of any
object, but also the two different orientations of $\mathbb R^d$ or
$\p^d$.

By the following theorem, the notions of chirotopes and hyperline
sequences are equivalent. This connects our concept
of hyperline sequences with other ways to look at oriented matroids.
The cyclic structure of a hyperline captures many instances of the
$3$--term Gra{\ss}mann--Pl\"ucker relations at once.  Therefore the
proof of the theorem becomes rather long and tedious.  But we
believe that this price is worth to pay.
Namely, it is easier to deal with a few cyclic structures than with a
multitude of Gra{\ss}mann--Pl\"ucker relations, and therefore it is
algorithmically more convenient and efficient to encode the structure of
oriented matroids in hyperline sequences rather than in chirotopes.

\begin{theorem}\label{thm:HLS-Chi}
  The set of positively oriented bases of a hyperline sequence of rank $r$ over
  $E$ is the set of positively oriented bases of a chirotope of rank $r$ over
  $E$, and vice versa.
\end{theorem}
\begin{proof}
\noindent 1. Let $X$ be a hyperline sequence of rank $r$ over $E$, and let
  $[x_1,\dots, x_r]$ be a positively oriented base of $X$.
  Since the cyclic orders
  $(x_{r-1},x_r,\overline{x_{r-1}},\overline{x_r})$ and 
  $(\overline{x_r}, x_{r-1},x_r,\overline{x_{r-1}})$ are equal up to a
  shift, $[x_1,\dots,x_{r-2},\overline{x_r},x_{r-1}]$ is a positively
  oriented base of $X$ as well. This together with Axiom~(H4) and an
  induction on the rank implies that all oriented simplices equivalent to
  $[x_1,\dots, x_r]$ are positively oriented bases of $X$, hence, yields
  Axiom~(C2).
  Axiom~(C1) follows from Axiom~(H1), Axiom~(C2) and induction on the
  rank.
  The Axioms~(H3) and~(C3) are equivalent. 
  
  In the next paragraphs we deduce Axiom~(C4) from the cyclic order of
  hyperlines. We use Axiom~(C2), that is already proven, but we do not
  mention any application explicitly.
  Let $x_1,x_2,\dots, x_r, y_1,y_2\in \mathbf E$
  satisfy the first two inequalities in Axiom~(C4).
  The third inequality is to prove. We can assume that $x_1,\dots,x_{r-2}$
  defines some hyperline $(Y|Z)\in X$ with period length $2k$, 
  since otherwise Axiom~(C4) is trivial.
  For simplicity, we write $[a,b]$ in place of $\chi([x_1,\dots,x_{r-2},
  a,b])$.
  
  If $[y_1,x_r]=[x_{r-1},y_2]=-1$ then we replace $y_1$ with $x_r$,
  $x_r$ with $y_1$, $x_{r-1}$ with $\overline{y_2}$ and $y_2$ with
  $\overline{x_{r-1}}$. This changes the signs of the factors in the first
  inequality, whereas the other inequalities remain unchanged.
  Similarly, we can assume without loss of generality that both 
  factors of the second inequality and $[y_1,y_2]$ are non-negative.
  Let $a,b,c,d \in C_{2k}$ so that $$x_{r-1}\in Z^a, x_{r}\in Z^b,
  y_1\in Z^c, y_2\in Z^d.$$
  By shifting the hyperline, we assume that
  $c=0$.  Then we have $$
  b, d-a, b-d, a, d \in \{0,1,\dots,k\}\subset
  C_{2k}.$$

  If $[y_1,y_2]=0$ then the third inequality is satisfied, and
  Axiom~(C4) is proven.  If $[y_1,y_2]=1$ then $d\in\{1,\dots,k-1\}$.
  Assume that $[x_{r-1},x_{r}]=-1$, thus $a-b\in \{1,\dots,k-1\}$.
  With the natural order on $0,1,\dots,k$ we find $0\le b < a\le d <
  k$. This is a contradiction to $b-d\in \{0,1,\dots,k\}$.  Hence
  $[x_{r-1},x_r]\ge 0$, which finishes the proof of Axiom~(C4).

\noindent
2.  Conversely, let $\chi$ be a chirotope of rank $r$ over $E$. We
wish to construct a hyperline sequence $X$ with the same positively
oriented bases. Note that in our proof we do not mention all
applications of Axiom~(C2) explicitly.  If $\chi$ is of rank $1$, then
we define $X$ as the set of all $x\in \mathbf E$ with $\chi([x]) =
+1$.  Since $\chi([x]) = -\chi([\overline x])$ and since
$X\not=\emptyset$ by Axiom~(C1), $X$ is a hyperline sequence of rank
$1$, and it has the desired positively oriented bases by construction.
  
Let $\chi$ be of rank $2$. 
Fix some element $e\in E$.  We iteratively
define a sequence $X^0,X^1,\ldots$ of subsets of $\mathbf E$ as
follows.
We start with
  \begin{eqnarray*}
    X^{0} = \big\{ x\in \mathbf E \;|\; &\hspace{-.5em}\chi([e,x])=1&
    \text{\hspace{-1em}, and 
    for all } y\in \mathbf E\\
      &&  \text{\hspace{-1em}with }\chi([e,y]) = 1\text{ holds } 
      \chi([x,y]) \ge 0\big\}. 
  \end{eqnarray*}
  For $a\ge 0$, pick some $x^{(a)}\in X^a$.  Define
  \begin{eqnarray*}
    X^{a+1}=
    \big\{ x\in \mathbf E \;|\; &\hspace{-.5em}\chi([x^{(a)},x])=1& \text{\hspace{-1em}, and
    for all } y\in \mathbf E\\
    && \text{\hspace{-1em}with }\chi([x^{(a)},y]) = 1\text{
    holds } \chi([x,y]) \ge 0\big\}. 
  \end{eqnarray*}
  We prove that $X^{a+1}$ (and similarly $X^0$) is not empty.  
  Under the assumption $X^{a+1}=\emptyset$, we inductively define 
  elements $x_0,x_1,\ldots\in \mathbf E$ as follows. 
  By Axiom~(C1), there is some $x_0\in \mathbf E$ with
  $\chi([x^{(a)},x_0])=1$.  For $i\ge 0$,
  since $x_i\not\in X^{a+1}$ there is some $x_{i+1}\in \mathbf E$ with
  $\chi([x^{(a)},x_{i+1}])=1$ and $\chi([x_i,x_{i+1}])=-1$.
  Since $\mathbf E$ is finite, we find an index $i\ge 2$ so that
  $x_{i+1}=x_k$ with $k<i-1$ (note that $\chi([x_i,x_{i-1}])=1$).
  We choose $i$  minimal, which implies $\chi([x_{i-1},x_k])\ge 0$.
  It follows
  \begin{eqnarray*}
    \chi([x^{(a)},x_k])\cdot \chi([x_i,x_{i-1}]) &=&1 \\
    \chi([x^{(a)},x_i])\cdot \chi([x_{i-1},x_k]) &\ge&0 \\
    \chi([x^{(a)},x_{i-1}])\cdot \chi([x_i,x_k]) &=&-1, 
  \end{eqnarray*}
  which contradicts Axiom~(C4).
  Hence $X^{a+1}\not=\emptyset$.
 Since $\chi([x^{(a)},x])=1$ for $x\in X^{a+1}$, we have
  $\chi([x^{(a)},\overline x])=-1$ by Axiom~(C2), thus $\overline
  x\not\in X^{a+1}$.
  In conclusion, $X^{a+1}$ is a hyperline sequence of rank $1$.

  We show that for $a>0$ and for all $x_1,y_2\in X^{a}$ and all
  $x_2\in \mathbf E$ we have $\chi([x_1,x_2]) = \chi([y_2,x_2])$.
  Let $y_1=x^{(a-1)}\in X^{a-1}$ be the element that appears in the
  definition of $X^a$. By definition
  $\chi([y_1,x_1])=\chi([y_1,y_2])=1$. Hence we have 
  $\chi([x_1,y_2])\ge 0$ and $\chi([y_2,x_1])\ge 0$, thus
  $\chi([x_1,y_2])=0$. 
  Assume that $\chi([x_1,x_2]) \not= \chi([y_2,x_2])$. Without loss of
  generality, we assume $\chi([x_1,x_2])=-1$ and $\chi([y_2,x_2])\ge
  0$.  We obtain
  \begin{eqnarray*}
    \chi([y_1,x_2])\cdot\chi([x_1,y_2]) &=& 0\\
    \chi([y_1,x_1])\cdot\chi([y_2,x_2]) &\ge& 0\\
    \chi([y_1,y_2])\cdot\chi([x_1,x_2]) &=& -1.
  \end{eqnarray*}
  This is a contradiction to Axiom~(C4), hence $\chi([x_1,x_2]) =
  \chi([y_2,x_2])$. This implies that $X^{a+1}$ does not depend on the
  choice of $x^{(a)}\in X^a$.

  Let $p\in \mathbb N$ be the least number so that $e\in X^{p-1}$.
  Since $X^{a+1}$ is independent of
  the choice of $x^{(a)}\in X^a$, it follows $X^a=X^{a+p}$ for all $a\in
  \mathbb N$.
  Moreover, if $q\in \mathbb N$ is the least number with $\overline
  e\in X^{q-1}$, we have $X^{a+q} = \overline{X^a}$ for all $a\in
  \mathbb N$. 

  In order to prove that $(X^0,X^1,\ldots, X^{p-1})$ yields a hyperline
  sequence of  rank 2 over $E$ with period length $2q=p$, it remains
  to show that $\mathbf E$ is a disjoint union of $X^0,\dots,X^{p-1}$.
  By choosing $p$ minimal and by the independence of $X^{a+1}$ from the
  choice of $x^{(a)}$, it follows that $X^0,\dots, X^p$ are disjoint.
  Let $x_1\in\mathbf E$. There is some $x_2\in \mathbf E$ with
  $\chi([x_1,x_2])=1$, by Axiom~(C1).
  Since $\chi([x^{(0)},x^{(1)}]) = 1$ by definition, Axiom~(C3) implies
  that $\chi([x_1,x^{(0)}])$ or $\chi([x_1,x^{(1)}])$ does not vanish.
  Therefore and by Axiom~(C2), there is a least index $i$ so that
  $\chi([x^{(i)},x_1])=1$ and $\chi([x^{(i+1)},x_1])\le 0$. 
  Let $y\in \mathbf E$ with $\chi([x^{(i)},y])=1$.
  By definition of $x^{(i+1)}$, we have 
  \begin{eqnarray*} 
    \chi([y,x^{(i)}])\cdot \chi([x^{(i+1)},x_1]) &\ge&0, \\
    \chi([y,x^{(i+1)}])\cdot \chi([x_1,x^{(i)}]) &\ge&0
  \end{eqnarray*}
  and $\chi([x^{(i)},x^{(i+1)}]) =1$, and therefore $\chi([x_1,y])\ge 0$
  by Axiom~(C4). But this means $x_1\in X^{(i+1)}$.
  In conclusion, $\mathbf E = X^0\cup\dots\cup X^p$.
  Hence $(X^0,X^1,\ldots, X^{p-1})$ yields a hyperline
  sequence  of  rank 2 over $E$, and by construction it has the same
  positively oriented bases than $\chi$. 
  
%%%%%%%%%%%%

  Finally, we come to the case of rank $r\ge 3$.  Let
  $x_1,\dots,x_{r-2} \in \mathbf E$ be so that there are two elements
  $x_{r-1},x_r\in \mathbf E$ with $\chi([x_1,\dots, x_r])=1$. Such an
  $(r-2)$--tuple exists by Axiom~(C1).
  It is easy to verify that 
  $$ C(x_1,\dots,x_{r-2}) = \{ [y_{r-1},y_r]\in \Delta_1(E) 
               |\; \chi([x_1,\dots, x_{r-2}, y_{r-1}, y_r])=1\}$$
  is the set of positively oriented bases of a chirotope of rank 2 over
  some subset of $E$, thus, of a hyperline sequence
  $Z(x_1,\dots,x_{r-2})$ of rank 2.
  By application of Axiom~(C2), one can also show that
  \begin{eqnarray*}
  B(x_1,\dots,x_{r-2}) = \{ [y_1,\dots, y_{r-2}]\in \Delta_{r-3}(E)
  &|& \chi([y_1,\dots, y_r])=1\\ &&\text{ for all } [y_{r-1},y_r]\in
  C(x_1,\dots,x_{r-2})\}
  \end{eqnarray*}
  is the set of positively oriented bases of
  a chirotope of rank $r-2$ over some subset of $E$,
  with $[x_1,\dots,x_{r-2}]\in B(x_1,\dots,x_{r-2})$.
  By induction, it is equivalent to a hyperline sequence 
  $Y(x_1,\dots,x_{r-2})$ of rank $r-2$.
  
  We collect all pairs $\big(Y(x_1,\dots,x_{r-2}) |
  Z(x_1,\dots,x_{r-2})\big)$ to form a set (not a multi-set)
  $X$, where $x_1,\dots,x_{r-2}\in \mathbf E$.  It has the same
  positively oriented bases as $\chi$, by construction. It remains to
  show that $X$ is a hyperline sequence of rank $r$ over $E$.
  The most difficult to prove is Axiom~(H1).
  Let $Y=Y(x_1,\dots,x_{r-2})$ and $Z=Z(x_1,\dots,x_{r-2})$. 
  By definition, in an oriented simplex over $E$ any element of $\mathbf
  E$ occurs at most once. This implies that $E(Y)\cap E(Z)=\emptyset$ .
  Let $e\in E\setminus E(Z)$. It remains to show that $e\in E(Y)$, hence
  $E=E(Y)\cup E(Z)$.  Let $[x_1,\dots,x_{r-2}]$ and $[x_{r-1},x_r]$ be
  positively oriented bases of $Y$ and $Z$, thus $\chi([x_1,\dots,
  x_r])=1$. By multiple application of Axioms~(C2) and~(C3), there is an index $i$
  so that $\chi([x_1,\dots,x_{i-1},e,x_{i+1},\dots,x_r])\not=
  0$. We have $i<r-1$ since $e\not\in E(Z)$.

  We claim that for all $[y_1,y_2]\in C(x_1,\dots,x_{r-2})$ holds
  $$\chi([x_1,\dots,x_{i-1},e,x_{i+1},\dots,x_{r-2},y_1,y_2]) =
   \chi([x_1,\dots,x_{i-1},e,x_{i+1},\dots,x_r]).$$
   By replacing $e$ with $\overline e$ if necessary, we can assume that
   $ \chi([x_1,\dots,x_{i-1},e,x_{i+1},\dots,x_r])=1$. 
   By Axiom~(C2), it suffices to consider the case $y_1=x_{r-1}$.  For
   simplicity, we denote $[a,b]$ for
   $\chi([x_1,\dots,x_{i-1},a,x_{i+1},\dots,x_{r-1},b])$.  We have
   $[x_i,x_r]=[x_i,y_2]=[e,x_r]=1$, and $[x_i,e]=0$
   since $e\not\in E(Z)$.
   If $[e,y_2]=0$, then
   \begin{eqnarray*}
     [x_r,y_2]\cdot [e,x_i] &=&0\\
     {}[x_r,x_i]\cdot [e,y_2] &=&0\\
     {}[x_r,e]\cdot [x_i,y_2] &=&-1, 
   \end{eqnarray*}
   in contradiction to Axiom~(C4). Hence $[e,y_2]\not=0$. Since
   \begin{eqnarray*}
     [x_i,y_2]\cdot [e,x_r] &=&1\\
     {}[x_i,e]\cdot [x_r,y_2] &=&0\\
   \end{eqnarray*}
   and $[x_i,x_r]=1$, it follows $[e,y_2]=1$ from Axiom~(C4).
   This proves our claim.
   The claim implies that $e\in E(Y)$, which finishes the proof of
   Axiom~(H1).

  We check the remaining axioms.  Since $\chi([x_1,\dots,x_r]) =
  \chi([x_1,\dots,x_{r-3},\overline{x_{r-2}}, x_r,x_{r-1}])$ by
  Axiom~(C2), we have
  \begin{eqnarray*}
    Y(x_1,\dots,x_{r-3},\overline{x_{r-2}}) &=& - Y(x_1,\dots,x_{r-2}),
    \\
    Z(x_1,\dots,x_{r-3},\overline{x_{r-2}}) &=& - Z(x_1,\dots,x_{r-2}),
  \end{eqnarray*}
  and this implies Axiom~(H2).
  Axiom~(H3) is Axiom~(C3). Finally, Axiom~(H4) is a special case of
  Axiom~(C2). Thus, $X$ is a hyperline sequence.
\end{proof}

Our proof of the topological representation theorem in
Sections~\ref{sec:FLbase} and~\ref{sec:FLgeneral} is by induction on
the rank and the number of elements of a chirotope or hyperline
sequence.  In the rest of this section, we expose the basic techniques
for this induction.
Let $\chi$ be a chirotope of rank $r$ over $E$, and let $R\subset E$.
We define the map $\chi\setminus R \co \Delta_{r-1}(E\setminus R)\to
\{-1,0,+1\}$ as the restriction of $\chi$ to $\Delta_{r-1}(E\setminus
R)$. We call $\chi\setminus R$ the \begriff{deletion} of $R$ in
$\chi$.  In general, $\chi\setminus R$ is not a chirotope. However if
the rank of $\chi$ is smaller than $|E|$, then we find by the
following lemma an element that can be deleted such that $\chi
\setminus \{i\}$ is a chirotope of the same rank.
\begin{lemma}\label{lem:fullrank}
  Let $\chi$ be a chirotope of rank $r$ over $E$.  If $|E|>r$ then
  there is an $i\in E$ so that $\chi\setminus \{i\}$ is a chirotope of
  rank $r$ over $E\setminus \{i\}$.
\end{lemma}
\begin{proof}
  By Axioms~(C1) and~(C2), there are $f_1,\dots, f_r \in E$ with
  $\chi([ f_1,\dots, f_r])\not=0$. Since $|E|>r$, we can pick some
  $i\in E\setminus \{f_1,\dots,f_r\}$.
  It is obvious that $\chi\setminus \{i\}$ satisfies Axioms~(C2), (C3)
  and~(C4).
  
  It remains to show that $\chi\setminus \{i\}$ satisfies Axiom~(C1).
  Let $e_1\in E\setminus \{i\}$. By Axiom~(C1), there are
  $e_2,\dots,e_r\in E$ with $\chi([ e_1,\dots, e_r])\not=0$. If
  $i\in\{e_2,\dots,e_r\}$ then by Axiom~(C2) we can assume that
  $i=e_r$.  Since $\chi([ f_1,\dots, f_r])\not=0$ from the preceding
  paragraph, we find some $f_j$ so that $\chi([ e_1,\dots, e_{r-1},
  f_j]) \not= 0$, by Axiom~(C3).
  Since $[ e_1,\dots, e_{r-1}, f_j] \in \Delta_{r-1}(E\setminus
  \{i\})$, the map $\chi\setminus\{i\}$ satisfies Axiom~(C1) and is
  therefore a chirotope.
\end{proof}

Let $R=\{e_1,\dots, e_k\}\subset E$, with $e_1<\dots < e_k$ and $k<r$. 
We define $E/R$ as the set of all $e \in E$ for which there exist
$e_{k+1},\dots e_{r}\in E$ so that $\chi([ e_1,\dots, 
e_{r-1},  e])\not= 0$. 
The map $\chi/R\co\Delta_{r-1-k}(E/R)\to \{-1,0,+1\}$ is then defined by 
$ \chi/R([ e_{k+1},\dots,  e_r]) = \chi([
e_1,\dots,  e_r])$, for all $[ e_{k+1},\dots,  e_r]
\in \Delta_{r-1}(E/R)$. 
It is obvious that $\chi/R$ satisfies Axioms~(C2), (C3) and~(C4). It
satisfies Axiom~(C1) by definition of $E/R$. Hence, $\chi/R$ is a
chirotope of rank $r-k$ over $E/R$, that is called the
\begriff{contraction} of $\chi$ on $R$. 

Finally, let $X$ be a hyperline sequence of rank $r$ over $E$. Let
$\chi$ be the chirotope that corresponds to $X$. For any $R\subset E$ so
that $\chi\setminus R$ (resp. $\chi/R$) is a chirotope, we define the
\begriff{deletion} $X\setminus R$ of $R$ in $X$ (resp. the
\begriff{contraction} $X/R$ of $X$ on $R$) as the hyperline sequence
associated to $\chi\setminus R$ (resp. $\chi/R$).

%******************************************************************************

\section{Arrangements of oriented pseudospheres}
\label{sec:PSA}

In the preceding sections, we extracted combinatorial data from vector
arrangements and turned properties of these data into axioms, yielding
hyperline sequences and chirotopes. In this section, we will generalize
vector arrangements in a geometric way. The main aim of this paper is to
prove the equivalence of these geometric structures with hyperline
sequences. 

The idea is as follows. Let $\p^{r-1}\subset \mathbb R^r$ be the unit
sphere, and let $V\subset R^r$ be a finite multiset of non-zero
vectors. Any $v\in V$ yields a pair $\pm \frac{v}{\|v\|}$ of points in
$\p^{r-1}$. This is dual to a hypersphere $S_v\in\p^{r-1}$, namely
the intersection of $\p^{r-1}$ with the hyperplane $H_v\subset \mathbb
R^r$ that is perpendicular to $v$. We put an
orientation on $H_v$, so that a positive base of $H_v$ together with $v$
is a positive base of $\mathbb R^r$. This induces an orientation on
$S_v$. 
In conclusion, $V$ is dual to an arrangement of oriented hyperspheres in
$\p^{r-1}$. The geometric idea is now to consider arrangements of
oriented embedded spheres of codimension $1$ in $\p^{r-1}$, intersecting
each other similarly to hyperspheres, though not being hyperspheres in
general. 

We formalize this idea.
Let $\p^d$ denote the $d$--dimensional oriented sphere 
$$\p^d = \left\{\left(x_1,\dots,x_{d+1}\right)\in\mathbb R^{d+1}\;|\;\;
                     x_1^2+\dots + x_{d+1}^2 = 1\right\},$$
and let
$$B^d=\big\{(x_1,\dots,x_{d}) \in \mathbb R^d | x_1^2+\dots +
 x_{d}^2 \le 1 \big\}$$
denote the closed $d$--dimensional ball.
A submanifold $N$ of codimension $m$ in a $d$--dimensional manifold $M$
is \begriff{tame} if any $x\in N$ has  an open neighbourhood
$U(x)\subset M$ such that there is a homeomorphism $\overline{U(x)}\to
B^d$ sending $U(x)\cap N$ to $B^{d-m}\subset B^d$.

An \begriff{oriented pseudosphere} $S\subset \p^d$ is a tame embedded 
$(d-1)$--dimensional sphere with a choice of an orientation. Obviously
any oriented hypersphere is an oriented pseudosphere. 
Let $\psi\co\p^{d-1}\to\p^d$ be an embedding with image $S$, inducing
the desired orientation of $S$. 
By a result of M.~Brown~\cite{brown2}, the image of $\psi$ is tame if
and only if $\psi$ can be extended to an orientation preserving
embedding  $$\tilde\psi\co \p^{d-1}\times [-1,1] \to \p^d \hbox{ with }
\psi(\cdot)=\tilde\psi(\cdot,0).$$
The image of an oriented pseudosphere $S$ under a homeomorphism
$\phi\co\p^d\to\p^d$ is again an oriented pseudosphere, 
since the defining embedding $\phi\circ\psi\co\p^{d-1}\to\p^d$ can be
extended to $\phi\circ\tilde\psi\co \p^{d-1}\times [-1,1] \to
\p^d$.
By the generalized Sch\"onflies theorem, that was also proven by
M.~Brown~\cite{brown}, $\p^d\setminus S$ is a disjoint union of two open
balls whose closures are closed balls.
We call the connected component of $\p^d\setminus S$ containing
$\tilde\psi(\p^{d-1}\times\{1\})$ (resp.\ $\tilde\psi(\p^{d-1}\times\{-1\})$)
the \begriff{positive side} $S^+$(resp.\ \begriff{negative side} $S^-$)
of $S$.
The following definition of  arrangements of oriented pseudospheres is
similar to~\cite{bjorner}, p.~227.  Recall $E_n=\{1,\dots, n\}$.

\begin{definition}
  Let $S_1,\dots, S_n\subset \p^d$ be not necessarily distinct oriented
  pseudospheres, ordered according to their indices.  
  Assume that the following conditions hold:
  \begin{itemize}
  \item[(A1)] $S_R= \p^d\cap\bigcap_{i\in R} S_i$ is empty or homeomorphic to a
    sphere, for all $R\subset E_n$. 
  \item[(A2)] Let $R\subset E_n$ and $i\in E_n$ with $S_R\not\subset
    S_i$. Then $S_R\cap S_i$ is a pseudosphere in $S_R$, and $S_R\cap  S_i^+$ 
    and $S_R\cap S_i^-$ are both non-empty. 
  \end{itemize}
  Then the ordered multiset $\{S_1,\dots, S_n\}$ is an \begriff{arrangement of
   oriented pseudospheres} over $E_n$. 
  The arrangement is called \begriff{of full rank} if the intersection
  of its members is empty.   
\end{definition}
We omitted Axiom~(A3) from~\cite{bjorner}, p.~227,  since it follows from the
other two axioms. We remark that in \cite{bjorner}, arrangements of full
rank are called \emph{essential}. Obviously if $S_1,\dots, S_n$ are
oriented hyperspheres then $\{S_1,\dots, S_n\}$ is an arrangement of
oriented pseudospheres. At the end of this section, we will characterize
arrangements of oriented pseudospheres by a single axiom.

\begin{definition}
  Two ordered multisets $\{S_1,\dots,S_n\}$ and $\{S'_1,\dots, S'_n\}$ of
  oriented pseu\-do\-sphe\-res in $\p^d$  are \begriff{equivalent} if
  there is an orientation preserving homeomorphism $\p^d\to\p^d$ sending
  $S_i^+$ to $(S_i')^+$ and $S_i^-$ to $(S_i')^-$, simultaneously 
  for all $i\in E_n$. 
  We do not allow renumbering of the pseudospheres. 
\end{definition}
Since the image of a pseudosphere under a homeomorphism $\phi\co \p^d \to
\p^d$ is a pseudosphere, it is easy to observe that if $\{S_1,\dots,
S_n\}$  is an arrangement of oriented pseudospheres, then so is
$\{\phi(S_1),\dots, \phi(S_n)\}$.

\begin{example}  
We construct non-equivalent arrangements $\mathcal A(d,+),
  \mathcal A(d,-)$ of $d+1$ oriented pseudo\-spheres
$S_1,\dots,S_{d+1}\subset \p^d$  of full rank as follows.
  For $i=1,\dots, d+1$, set
  $$ S_i = \left\{ (x_1,\dots,x_{d+1})\in \p^d\;|\;\; x_i=  0\right\}.$$
  In   $\mathcal A(d,+)$, define $S_1^+ = 
  \big\{ (x_1,\dots,x_{d+1})\in \p^d\;|\;\; x_1>  0\big\}$. In
  $\mathcal A(d,-)$, define $S_1^+ = \big\{ (x_1,\dots,x_{d+1})\in
    \p^d\;|\;\; x_1 <  0\big\}$. In both $\mathcal A(d,+)$
  and$\mathcal A(d,-)$,  define  
  $S_i^+ =\big\{ (x_1,\dots,x_{d+1})\in \p^d\;|\;\; x_i >  0\big\}$, for
  $i=2,\dots,d+1$. 

  $\mathcal A(d,+)$ and $ \mathcal A(d,-)$ are of full
  rank. If $\phi\co\p^d\to\p^d$ is any homeomorphism that fixes
  $S_2^+,\dots,S_{d+1}^+$ setwise and maps  $\big\{
    (x_1,\dots,x_{d+1})\in \p^d\;|\;\; x_1>  0\big\}$ to $\big\{
    (x_1,\dots,x_{d+1})\in \p^d\;|\;\; x_1<  0\big\}$, then $\phi$ is
  orientation reversing. Thus, $\mathcal A(d,+)$ is not equivalent to
  $\mathcal A(d,-)$. 
\end{example}

In the preceding section, we have defined deletion and contraction of
chirotopes and hyperline sequences. There is a similar notion for
arrangements of oriented pseudospheres, as follows.
Fix an arrangement $\mathcal A = \{S_1,\dots S_n\}$  of oriented
pseudospheres in $\p^d$. 
For any $R \subset E_n$, we obtain an arrangement $$\mathcal
A\setminus \big\{S_{r} |\; r\in R\big\}$$ of oriented pseudospheres over
$E_n\setminus R$ in $\p^d$ (Axioms (A1) and (A2) are easy to verify). We
denote this arrangement by $\mathcal A\setminus R$ and call 
it the \begriff{deletion} of $R$ in $\mathcal A$. In general,
$\mathcal A\setminus R$ is not of full rank, even if $\mathcal A$ is.

It is intuitively clear from Axiom~(A2), that for $R\subset E_n$ one
gets an arrangement of oriented pseudospheres on $S_R$, induced by
$\mathcal A$. 
In the following iterative definition of this induced arrangement, the
orientation of $S_R$ requires some care.
Let $r\in E_n$, and let $\psi_r\co \p^{d-1}\to\p^d$ be a tame embedding
defining $S_r$ with the correct orientation.
Denote $S_i' = \psi_r^{-1}(S_i\cap S_r)$, for $i\in E_n$.  We obtain
an ordered multiset
$$ \mathcal A/\{r\} = \{ S_i' |\; i\in E_n, S_r\not\subset S_i \}$$
of oriented pseudospheres in $\p^{d-1}$, where $(S'_i)^+ =
\psi_r^{-1}(S^+_i\cap S_r)$. Axioms~(A1) and~(A2) are easy to verify,
thus, $\mathcal A/\{r\}$ is an arrangement of oriented pseudospheres
over a subset of $E_n\setminus \{r\}$.

This construction can be iterated. Let $R\subset E_n$ so that $\dim S_R
= d - |R|$. List the elements of $R$ in ascending order, $r_1< r_2< \dots<
r_{|R|}$. 
The \begriff{contraction} of $\mathcal A$ in $R$ is then the arrangement
of oriented pseudospheres 
$$\mathcal A/R = \Big(\dots\big((\mathcal A/\{r_1\}) \big/
\{r_2\}\big)\dots \Big/ \{r_{|R|}\}\Big).$$
Note that $(\mathcal A/\{r_1\})/\{r_2\}$ and $(\mathcal
A/\{r_2\})/\{r_1\}$ are related by an orientation \emph{reversing}
homeomorphism $\p^{d-2}\to\p^{d-2}$.
It is easy to see that $\mathcal A/ R$ is of full rank if and only if
$\mathcal A$ is of full rank. 

An arrangement $\mathcal A=\{S_1,\dots, S_n\}$ of oriented pseudospheres
in $\p^d$ yields a cellular decomposition of $\p^d$, as follows.
For any subset $I\subset E_n$, we consider the parts of the intersection
of the pseudospheres with label in $I$ that are not contained in
pseudospheres with other labels,
\begin{eqnarray*}
 \C(I,\mathcal A) &=& \big\{x \in \p^d | \; x\in S_i \iff i\in I \text{, for
  all } i \in E_n\big\}\\
  &=& \big(\bigcap_{i\in I} S_i\big)\setminus \big(\bigcup_{j\in
    E_n\setminus I} S_j\big).
\end{eqnarray*} 
By the next theorem, the connected
components of $\C(I,\mathcal A)$ are topological cells if $\mathcal A$
is of full rank. 
Hence, it makes sense to refer to a connected component of
$\C(I,\mathcal A)$ as a \begriff{cell of $\mathcal A$ with label $I$}.   The 
\begriff{$d$--skeleton} $\mathcal A^{(d)}$ of $\mathcal A$ is the union of 
its cells of dimension $\le d$. 
By Axiom (A1), $\bigcap_{i\in R} S_i$ is empty or a sphere, for
$R\subset E_n$. Thus, for any
$0$-dimen\-sio\-nal cell of $\mathcal A$  exists exactly one other
cell of $\mathcal A$ with the same label, corresponding to the two
points in $\p^0$.  
The statement that the cells of an arrangement of oriented pseudospheres
of full rank are in fact topological cells is essential for the proof of
the Topological Representation Theorem in Sections~\ref{sec:FLbase}
and~\ref{sec:FLgeneral}.  

\begin{theorem}\label{thm:cells}
  Let $\mathcal A\not=\emptyset$ be an arrangement of oriented
  pseudospheres over $E_n$. Any connected component of
  $\C(\emptyset,\mathcal A)$  is a $d$--dimensional cell whose closure
  is a closed ball.  If $\mathcal A$ is of full rank, then for any
  $I\subset E_n$, any connected component of  $\C(I,\mathcal A)$ is an
  open cell whose closure is a closed ball.
\end{theorem}
\begin{proof}
  First, we prove that the theorem holds for $\C(\emptyset,\mathcal A)$,
  by induction on the number $n$ of elements of $\mathcal A$. 
  Since $\mathcal A \not=\emptyset$, we have $n>0$. 
  The base case $n=1$ is the generalized Sch\"onflies
  theorem~\cite{brown}, stating that an embedded
  tame $(d-1)$--sphere in $\p^d$ is the image of a hypersphere under a
  homeomorphism $\p^d\to\p^d$. 
  In particular, the complement of a pseudosphere is a disjoint union of
  two $d$--dimensional cells whose closures are balls.

  Let $n>1$, and let $c$ be the closure of a
  connected component of $\C(\emptyset,\mathcal A\setminus \{n\})$. By
  induction on $n$, it is a cell of dimension $d$.
  The  connected components of $S_n\cap c$ correspond to the closures of
  connected components of $\C(\emptyset, \mathcal A/\{n\})$. By
  induction hypothesis, applied to the arrangement $\mathcal A/\{n\}$ of 
  $n-1$ pseudospheres in $\p^{d-1}$, these are closed balls of dimension
  $d-1$ whose boundary lies in $\d  c$.  
  Since $S_n$ is tame, $S_n\cap c$ is tame in $c$. Hence, it follows
  from the generalized Sch\"onflies theorem~\cite{brown} that the
  closure of any connected component of $c\setminus S_n$ is a ball of
  dimension $d$.
% [more details!].  
  In conclusion,
  the connected components of $\C(\emptyset,\mathcal A)$ are cells of
  dimension $d$ whose closures are balls.

  It remains to prove the theorem in the full rank case with
  $I\not=\emptyset$. 
  We can assume $\C(I,\mathcal A)\not=\emptyset$.
  Since $\mathcal A$ is of full rank, $\C(E_n, \mathcal
  A)=\emptyset$, hence $I\not= E_n$.
  By Axiom~(A1), $S_I$ is homeomorphic to some
  sphere $\p^e$. There is an $R\subset I$ with $|R|\le d-e$ and
  $S_R = S_I$. 
  The set $\C(I,\mathcal A)$  is mapped to $\C(\emptyset,  \mathcal A/R)$
  by the restriction of a homeomorphism $S_I\to \p^e$. 
  It has already been proven that the closure of any connected component  of
  $\C(\emptyset,  \mathcal A/R)$ is a ball. Thus, the closure of any connected
  component of  $\C(I,\mathcal A)$ is a ball, as well. 
\end{proof}

We remark that the preceding theorem becomes wrong by dropping the
hypothesis that pseudospheres are tame. In fact, there are wild
$2$-spheres in $\p^3$ (e.g. the famous Horned Sphere of Alexander),
whose complement is not a union of cells. We will consider arrangements
of oriented embedded spheres (not necessarily tame) in
Section~\ref{sec:wild} and prove that these wild arrangements have the
same combinatorics than tame arrangements.

Let $\mathcal A = \{S_1,\dots, S_n\}$ be an ordered multiset of oriented
pseudospheres in $\p^d$. For $R\subset E_n$, denote $\mathcal A_R =
\{S_j | j\in R\}$.
In the remainder of this section, we show that one can replace
Axioms~(A1) and~(A2) by the following single \begriff{Axiom~(A')}. 
\begin{itemize}
\item[(A')] Let $R\subset E_n$ so that $S_{R'}\not= S_R$ for any proper
  subset $R'$ of $R$.
  Then, $\mathcal A_{R}$ is equivalent to an arrangement of
  $|R| $ oriented hyperspheres in $\p^d$.
  \end{itemize}

\begin{theorem}\label{thm:PSA}
  An ordered multiset $\mathcal A= \{S_1,\dots, S_n\}$ of oriented
  pseudospheres in $\p^d$ satisfies Axiom~(A') if and only if it is an
  arrangement of oriented pseudospheres.
\end{theorem}
\begin{proof}
  First, we assume that $\mathcal A$ satisfies Axiom~(A') and prove that
  $\mathcal A$ is an arrangement of oriented pseudospheres. 
  Let $\emptyset \not= R\subset E_n$, and assume that $S_R\not=\emptyset$. 
  Replacing $R$ by a subset if necessary, we can assume
  that $S_{R'}\not= S_R$ for any proper subset $R'$ of $R$. 
  By Axiom~(A'),
  $\mathcal A_R$ is equivalent  to an arrangement of oriented
  hyperspheres. Thus, $S_R$ is homeomorphic to a sphere, and Axiom~(A1)
  holds. 
  Let $i\in E_n$ with $S_R\not\subset S_i$; we wish to prove (A2).
  There is some $\tilde R\subset R\cup\{i\}$ with $S_{\tilde
    R}=S_{R\cup\{i\}}$, so that $S_{R'}\not=S_{\tilde R}$ for any proper
  subset $R'$ of $\tilde R$.
  Since by Axiom~(A'), $\mathcal A_{\tilde R}$ is equivalent to an
  arrangement of oriented hyperspheres, it follows
  $$ \dim S_{\tilde R\setminus \{i\}} = \dim S_{R\cup \{i\}} +1.$$
  Since both $S_R$ and $S_{R\cup \{i\}}$ are spheres by Axiom~(A1) (that
  has already been proven) and since  $S_R\not\subset S_i$ by hypothesis, it
  follows 
  $$
  \dim S_{R\cup \{i\}} +1 = \dim S_{\tilde R\setminus \{i\}} \ge \dim
  S_R > \dim S_{R\cup\{i\}}.$$
  Thus, $S_{\tilde R\setminus \{i\}}\subset S_R$ is a pair of speres of
  the same dimension $\dim S_{R\cup \{i\}} +1$, hence, $S_{\tilde
    R\setminus \{i\}}= S_R$.  
  Since $\mathcal A_{\tilde R}$ is equivalent to an arrangement of
  oriented hyperspheres, Axiom~(A2) holds for $S_{\tilde R\setminus
    \{i\}}= S_R$ and $S_i$.  In conclusion, $\mathcal A$ is an
  arrangement of oriented pseudospheres.
  
  Secondly, we assume that $\mathcal A$ is an arrangement of oriented
  pseudospheres and prove by induction on $|R|$ that (A') holds. If
  $|R|=1$, then (A') is nothing but 
  the generalized Sch{\"o}nflies  theorem~\cite{brown}.
  In the general case, let $\emptyset \not= R\subset E_n$ so that $S_{R'}\not= S_R$ for
  any proper subset $R'$ of $R$, and let $i\in R$.
  It follows easily that $S_{R'}\not=S_{R\setminus \{i\}}$ for any
  proper subset $R'$ of $R\setminus \{i\}$. Thus, by Axiom~(A') we can
  assume that $\mathcal A_{R\setminus\{i\}}$ is an arrangement of
  oriented hyperspheres. Note that the connected components of $\mathcal
  C(I,\mathcal A_{R\setminus \{i\}})$ are topological cells, for any
  proper subset $I$ of $R\setminus \{i\}$.
  
   Let $H\subset \p^d$ be a hypersphere that does not contain 
  $S_{R\setminus\{i\}}$. Our aim is to transform $S_i$ into $H$, fixing
  $\mathcal A_{R\setminus\{i\}}$ cellwise, which implies Axiom~(A') for
  $\mathcal A_R$.
  By Axiom~(A2) and by the generalized Sch\"onflies theorem,
  $S_{R\setminus \{i\}}\cap
  S_i$ can be mapped to $S_{R\setminus\{i\}}\cap H$ by some orientation
  preserving homeomorphism $S_{R\setminus\{i\}}\to S_{R\setminus \{i\}}$,
  The homeomorphism can be extended to a homeomorphism $\p^d\to\p^d$
  fixing all cells of $\mathcal A_{R\setminus\{i\}}$, by
  the cone construction~\cite{stoeckerzieschang}. 

  We now proceed with transforming $S_i$ in cells of $\mathcal
  A_{R\setminus\{i\}}$ of  higher dimension. 
  Let $R'$ be a proper subset of $R\setminus \{i\}$.  By induction on
  $|R\setminus R'|$, we can assume that 
  $S_{R'\cup \{j\}}\cap S_i = S_{R'\cup\{j\}}\cap H$, for all $j\in
  R\setminus (R'\cup \{i\})$.
  Let $C\subset S_{R'}$ be the closure of a cell of $\mathcal
  A_{R\setminus \{i\}}$ with $\dim C = \dim S_{R'}$.  It follows from
  Axiom~(A2) that $B=S_i\cap C$ is a tame ball in $C$. 
  As a consequence of the generalized Sch\"onflies
  theorem and since $\d B\subset H$, we can map $C\cap S_i$ to $C\cap H$
  by an orientation preserving homeomorphism $C\to C$ that is the
  identity on $\d C$. It can be extended to a homeomorphism $\p^d\to
  \p^d$ by the cone construction, fixing all cells of $\mathcal
  A_{R\setminus \{i\}}$. 
  Thus, we can transform $S_{R'}\cap S_i$ into $S_{R'}\cap 
  H$. 
  Finally, when we achieve $R'=\emptyset$, we have $S_{R'}\cap S_i = S_i
  = H$.
  Therefore, $\mathcal A_{R}$ is equivalent to an arrangement
  of oriented hyperspheres. This proves (A'), as claimed.
\end{proof}

\begin{corollary}\label{cor:n=d+1}
  If an arrangement $\mathcal A$ of oriented pseudospheres in $\p^{d}$ is
  of full rank and $n=d+1$, then $\mathcal A$ 
  is equivalent to $\mathcal A(d,+)$ or $\mathcal A(d,-)$.
\end{corollary}
\begin{proof}
  Let $R\subset E_{d+1}$ be minimal so that $\mathcal A_R$ is of full
  rank. 
  By the preceding theorem, $\mathcal A_R$ is equivalent to
  an arrangement  of oriented hyperspheres of full rank. Such an
  arrangement consists of at least $d+1$ hyperspheres, hence
  $R=E_{d+1}$ and $\mathcal A_R=\mathcal A$.

  Any  arrangement of $d+1$ oriented hyperspheres of full rank is dual to an
  arrangement of $d+1$ unit vectors  in $\mathbb R^{d+1}$ that span
  $\mathbb R^{d+1}$. 
  Since $GL_{d+1}(\mathbb R)$ acts transitively on those vector
  arrangements, it follows that $\mathcal A$ is equivalent to  $\mathcal
  A(d,+)$ or $\mathcal A(d,-)$, depending on its orientation.  
\end{proof}

%******************************************************************************

\section{Chirotopes and hyperline sequences associated to arrangements
  of oriented pseudospheres} 
\label{sec:X(A)}

The aim of our paper is to prove a topological representation theorem
for hyperline sequences, i.e., to establish a one-to-one correspondence
between hyperline sequences and equivalence classes of arrangements of
oriented pseudospheres.
In this  section, we settle one direction of this correspondence. We
associate to any arrangement of oriented pseudospheres a 
hyperline sequence, compatible with deletions and contractions.

We first expose the geometric idea.
By a \begriff{cycle} of an arrangement $\mathcal A$ of oriented
pseudo\-spheres in $\p^d$, we mean an embedded  circle 
$\mathbb S^1\subset \p^d$ that is the intersection of some
elements of $\mathcal A$. 
Let $L$ be a cycle of $\mathcal A$ with a choice of an orientation. It
corresponds to a hyperline $(Y|Z)$, as follows. The positively oriented
bases of $Y$ correspond to 
$(d-1)$--tuples of pseudospheres containing $L$ so that the system of
positive normal vectors of the pseudospheres together with a positive
tangential vector of $L$ forms a direct base of the oriented vector
space $\mathbb R^d$. 
The set of all elements of $\mathcal A$ containing $L$
corresponds to $E(Y)$. 
The $0$-dimensional cells of $\mathcal A$ on $L$
occur in a cyclic order, corresponding to the cyclic order of $Z$. 
Let $S_e\in \mathcal A$. If in a point of $L\cap S_e$ the
cycle $L$ passes from the negative side of $S_e$ to 
\begin{figure}[htbp]
  \begin{center}
    \epsfig{file=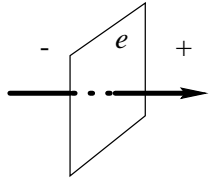}
    \caption{The atom $\{ e\}$ on a cycle}
    \label{fig:positivevertex}
  \end{center}
\end{figure}
the positive side, then we have an element $ e$ in the corresponding
atom of $(Y|Z)$ (see Figure~\ref{fig:positivevertex}). Further, in
the second point of $L\cap S_e$, the cycle $L$ passes from the positive
to the negative side,  yielding an element $\overline{ e}$ in the
atom that is opposite to the first atom.  

The meaning of Axiom~(H3) in the setting of arrangements of oriented
pseudospheres is that any two cycles of $\mathcal A$ have non-empty
intersection.
\begin{figure}[htbp]
  \begin{center}
    \epsfig{file=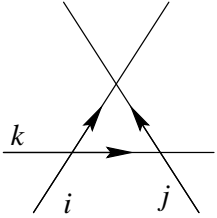}
    \caption{Axiom (H4) in rank 3}
    \label{fig:G0}
  \end{center}
\end{figure}
Figure~\ref{fig:G0} provides a visualization of Axiom~(H4) in an
arrangement of three oriented pseudospheres in $\p^2$:  If we get from
the pseudosphere $ k$ the cyclic order $\big(\{ i\}, \{ j\},
\{\overline{\imath}\}, \{\overline{\jmath}\}\big)$, then we get the cyclic order
$\big(\{\overline{k}\}, \{ j\}, \{ k\},\{\overline{\jmath}\}\big)$ from $i$
and the cyclic order $\big(\{\overline{ k}\}, \{\overline{\imath}\},
\{k\},\{i\}\big)$ from $ j$.

We formalize this idea in the rest of this section.
Let $\mathcal A=\{S_1,\dots,S_n\}$ be an arrangement of $n$ oriented
pseudospheres of full rank in $\p^d$. Our aim is to associate to
$\mathcal A$ a hyperline sequence $X(\mathcal A)$ of rank $d+1$ over
$E_n$. By the equivalence of hyperline sequences and chirotopes,
established in Theorem~\ref{thm:HLS-Chi}, this also allows to associate a
chirotope $\chi(\mathcal A)$ to $\mathcal A$, with the same
positively oriented bases.
If $d=0$ then define
$$ X(\mathcal A) = \big\{e\in E_n|\; S_e^+ = \{+1\}\big\} 
                           \cup \big\{\overline e\in \overline{E_n} |\;
                           S_e^+ = \{-1\}\big\},$$
which is obviously a hyperline sequence over $E_n$ of rank 1.

In the case $d=1$, the orientation of $\p^1$ yields a cyclic order
$p_0,p_1,\dots, p_{2k-1}$ on the points of $S_1\cup\dots\cup S_n$.  For
$a\in \{0,\dots,2k-1\}$, we define $X^a\subset \mathbf{E_n}$ by
\begin{enumerate}
\item $e\in X^a$ if $p_a\in S_e$ and, along the cyclic orientation of
  $\p^1$, one passes in $p_a$ from $S_e^-$ to $S_e^+$, and 
\item $\overline e\in X^a$ if $p_a\in S_e$ and  one passes in $p_a$ from
  $S_e^+$ to $S_e^-$. 
\end{enumerate}
It is easy to check that $(X^0,\dots,X^{2k-1})$ yields a hyperline
sequence of rank 2 over $E_n$.

In the case $d=2$, let $\gamma\subset \p^2$ be an oriented cycle of
$\mathcal A$. Let $R\subset E_n$ be the indices of oriented
pseudospheres containing $\gamma$. 
There is a hyperline sequence $Y(\gamma)$ of rank 1 over $R$, with $e\in
Y(\gamma)$ (resp. $\overline e\in Y(\gamma)$) if the orientation of
$S_e$ coincides (resp. does not coincide) with the orientation of
$\gamma$. 
As in the preceding paragraph, we
obtain a hyperline sequence $Z(\gamma)$ of rank 2 over $E_n\setminus
R$. 
We collect the pairs $(Y(\gamma),Z(\gamma))$ to form a set $X(\mathcal
A)$, where $\gamma$ runs over all oriented cycles of $\mathcal A$. 
Axioms~(H1) and~(H2) are obvious for $X(\mathcal A)$. To prove
Axiom~(H3), let $[x_1,x_2,x_3]$ and $[y_1,y_2,y_3]$ be two positively
oriented bases of $X(\mathcal A)$. By definition, the pseudospheres
$S_{y^*_1}$, $S_{y^*_2}$, $S_{y^*_3}$ have no point in common. In
particular, one of them, say,  $S_{y^*_1}$, does not contain
$S_{x^*_1}\cap S_{x^*_2}$. Thus $S_{y^*_1}$ intersects $S_{x^*_1}$
transversely in $S_{x^*_1}\setminus S_{x^*_2}$, hence $[x_1,x_2,y_1]$ or
$[x_1,x_2,\overline {y_1}]$ is a positively oriented base of $X(\mathcal
A)$. 
It remains to prove Axiom~(H4). Here we use the Jordan--Sch\"onflies
theorem~\cite{stoeckerzieschang}, stating that the complement of an
embedded $1$--sphere in $\p^2$ is a disjoint union of two discs. This
holds even without the assumption of tameness.
With this in mind, Axiom~(H4)  can be read off from Figure~\ref{fig:G0}.

In the case $d\ge 3$, let $\gamma$ be an oriented  cycle of $\mathcal
A$.
Let $R_\gamma= \{r\in E_n |\; \gamma\subset S_r\}$. 
As in the preceding paragraph, the cyclic orientation of $\gamma$
induces a cyclic order of the oriented points 
$$\gamma\cap \bigcup_{e\in E_n\setminus R_\gamma} S_e,$$ 
yielding a hyperline sequence $Z(\gamma)$ of rank 2 over $E_n\setminus
R_\gamma$. 
Since $\mathcal A$ is of full rank, there are $i,j\in E_n$ so that
$S_i\cap S_j$ is a sphere of dimension $d-2$ disjoint to $\gamma$.  We
may assume that $Z(\gamma)$ yields the cyclic order
$(\{i\},\{j\},\{\overline \imath\},\{\overline \jmath\})$, by changing the roles
of $i$ and $j$ if necessary.
Then $$\mathcal A(i,j) = \big((\mathcal
A/\{i\})\big/\{j\}\big)_{R_\gamma}$$
is an arrangement of oriented
pseudospheres in $\p^{d-2}$ of full rank over $R_\gamma$, and by
induction it corresponds to a hyperline sequence $Y(i,j)$ of rank
$d-1$ over $R_\gamma$.

We show that $Y(i,j)$ does not depend on the choice of $i,j$.  By
symmetry of $i$ and $j$, it suffices to pick $k\in E_n$ so that $S_i\cap
S_k$ is a sphere of dimension $d-2$ disjoint to $\gamma$ and $Z(\gamma)$
yields the cyclic order $(\{i\},\{k\},\{\overline\imath\},\{\overline
k\})$, and to show that $Y(i,j)=Y(i,k)$.
We will prove that the cyclic order of signed points on oriented cycles
of $\mathcal A(i,j)$ coincides with those of $\mathcal A(i,k)$.  It is
clear that this implies $Y(i,j)=Y(i,k)$.

\begin{figure}[htbp]
  \begin{center}
    \epsfig{file=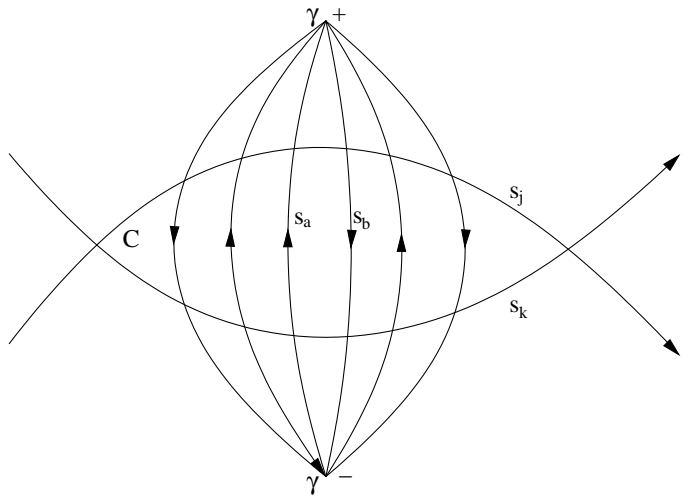, width=9cm}
    \caption{Parallel cycles}
    \label{fig:parallel}
  \end{center}
\end{figure}

Let $R\subset R_\gamma$ so that $S_R\cap S_i\cap S_j$ is a cycle of
$\mathcal A(i,j)$. We consider the $2$--sphere $S=S_R\cap S_i$. Both
$s_j=S\cap S_j$ and $s_k=S\cap S_k$ are embedded $1$--spheres in $S$
that are either equal or intersect in two points. Both $s_j$ and $s_k$
have the positive (resp. negative) point of $\gamma\cap S$ on their
positive (resp. negative) side. Thus, again using the
Jordan--Sch\"onflies theorem, the situation is as in Figure~\ref{fig:parallel}.
Let $C$ be a connected component of $S\setminus (s_j\cup s_k)$ disjoint
from $\gamma$. Let $a,b\in R_\gamma\setminus R$.  
Since $s_a=S\cap S_a$ is equal to $s_b=S\cap S_b$ or $s_a$ intersects
$s_b$ transversaly in $S\cap \gamma$, it follows that $C\cap s_a$ and
$C\cap s_b$ are parallel arcs. Hence, the cyclic order of the four
signed points of $s_j\cap (s_a\cup s_b)$ coincides with the order on
$s_k\cap (s_a\cup s_b)$. This proves our claim $Y(i,j)=Y(i,k)$.
Since  $Y(i,j)$ does not depend on the choice of $i$ and $j$, we denote
$Y(\gamma)=Y(i,j)$. 

Now, we define $X(\mathcal A)$ as the set of all pairs
$\big(Y(\gamma)|Z(\gamma)\big)$, where $\gamma$ runs over all oriented
cycles of $\mathcal A$, each cycle occuring in both orientations.
It remains to show that $X(\mathcal A)$ is a hyperline sequence of rank
$d+1$ over $E_n$. 
First of all, $X(\mathcal A)$ is not empty. It remains to check the four
axioms.
Axiom~(H1) is trivial.
Any sub-arrangement of $d-1$ pseudospheres defines a cycle $\gamma$ of
$\mathcal A$, up to orientation. Hence Axiom~(H2) follows from
$Y(\gamma)= - Y(-\gamma)$ and $Z(\gamma)= -Z(-\gamma)$.
Axioms~(H3) and~(H4) follow from the corresponding axioms in the case
$d=2$, since we can argue by contraction onto a $2$--sphere containing
the two cycles of $\mathcal A$ involved in Axioms~(H3) and~(H4).

In conclusion, $X(\mathcal A)$ is a hyperline sequence. By
Theorem~\ref{thm:HLS-Chi}, we can define $\chi(\mathcal A)$ as the
chirotope that is associated to the hyperline sequence $X(\mathcal
A)$ with the same positively oriented bases.  Let $R\subset E_n$ so
that $\mathcal A\setminus R$ (resp. $\mathcal A/R$) is an arrangement
of oriented pseudospheres of full rank. It is easy to see that
$\chi(\mathcal A\setminus R) = \chi(\mathcal A)\setminus R$ (resp.
$\chi(\mathcal A/R) = \chi(\mathcal A)/R$).

%********************************************************************************

\section{The topological representation theorem --- statement and base case}
\label{sec:FLbase}

\begin{theorem}[Topological Representation Theorem]\label{thm:FL}
  To any hyperline sequence $X$ of rank $r$ over $E_n$, there is
  an arrangement $\mathcal A(X)$ of  $n$ oriented pseudo 
  hyperspheres in $\mathbb S^{r-1}$ of full rank with $X=X(\mathcal A(X))$. The
  equivalence class of $\mathcal A(X)$ is unique. 
\end{theorem}

We prove Theorem~\ref{thm:FL} by induction on the number of elements and
the rank of $X$.  
In this section, we prove the base cases $r\le 2$ and $n=r$. The next
section is devoted to the inductive step.

\begin{lemma}\label{lem:FLbasecase_r=1}
  Theorem~\ref{thm:FL} holds for $r=1$.
\end{lemma}
\begin{proof}
  Let $X$ be a hyperline sequence of rank $1$ over $E_n$. We will
  construct an arrangement $\mathcal A(X)$ of oriented pseudospheres in
  $\p^0$, with $X(\mathcal A(X)) = X$. 

  An oriented pseudosphere $S$ in $\p^0= \{+1,-1\}$ is the empty set,
  together with the information whether $S^+=\{+1\}$ or $S^+=\{-1\}$. 
  Let $\mathcal A(X)= \{S_1,\dots, S_n\}$ be defined as follows.
  For $i\in E_n$, set $S_i^+= \{+1\}$ if $ i\in X$ and $S_i^+=
  \{-1\}$ otherwise. It is obvious that $X(\mathcal A(X)) = X$ and that 
  $\mathcal A(X)$ is unique with this property.
\end{proof}

\begin{lemma}\label{lem:FLbasecase_r=2}
  Theorem~\ref{thm:FL} holds for $r=2$.
\end{lemma}
\begin{proof}
  Let $X$ be a hyperline sequence of rank $2$ over $E_n$. 
  We will construct an arrangement $\mathcal A(X)$ of oriented
  pseudospheres in $\p^1$ with $X(\mathcal A(X)) = X$, and show that it
  is unique up to equivalence.

  The hyperline sequence  $X$ is a map from some cyclic group $C_{2k}$
  to non-empty subsets of $\mathbf E_n$. We consider $C_{2k}$ as a
  subgroup of $\p^1$.  
  An oriented pseudosphere corresponds to an embedding of
  $\p^0=\{+1,-1\}$ into $\p^1$. 
  For $i\in E_n$, let $a_i\in C_{2k}$ so that $ i\in X^a$. Define
  $\psi_i(+1) = a_i\in \p^1$ and $\psi_i(-1) = -a_i \in \p^1$. 
  It follows easily that $X(\mathcal A(X)) = X$ and that 
  $\mathcal A(X)$ is essentially unique with this property.
\end{proof}

\begin{lemma}\label{lem:n=r}
  There are exactly two chirotopes of  rank $|E|$ over $E$, namely one
  with $[ 1,\dots,  r]$ as positively oriented base, and
  the other with $[ 1,\dots,  {r-1}, \overline{ r}]$ as
  positively oriented base. 
\end{lemma}
\begin{proof}
  Set $r= |E|$. Without loss of generality, let $\chi$ be a chirotope of
  rank $r$ over $E_r=E$.
  There are exactly two equivalence classes of oriented
  $(r-1)$--simplices in $E_r$, namely those equivalent to $[ 1,\dots,
   r]$ and those equivalent to $[ 1,\dots, {r-1},
  \overline{ r}]$. 
  Thus, by Axiom~(C2), $\chi$ is completely determined by $\chi([
  1,\dots, r])$. 
  Axiom~(C1) implies  $\chi([ 1,\dots,  r])=\pm 1$. Hence there are at
  most two chirotopes of rank $|E|$ over $E$. 
  The conditions in Axioms~(C3) and~(C4) are empty for $|E|=r$. Thus,
  there are two chirotopes of rank $|E|$ over $E$. 
\end{proof}

\begin{lemma}\label{lem:FLbasecase_n=r}
  Theorem~\ref{thm:FL} holds for $n=r$.
\end{lemma}
\begin{proof}
  We prove the lemma for chirotopes rather than hyperline sequences.
  Let $\chi$ be a chirotope of rank $r$ over $E_r$.
  By Lemma~\ref{lem:n=r}, $\chi$ is determined by whether
  $[ 1,\dots, r]$ is a positively oriented base of $\chi$
  or not. In the former case, define $\mathcal A(\chi) = \mathcal
  A(r-1,+)$, in the latter case define $\mathcal A(\chi)= \mathcal
  A(r-1,-)$.
  We have $\chi(\mathcal A(\chi)) = \chi$ by construction, and the
  uniqueness of $\mathcal A(\chi)$ follows from Corollary~\ref{cor:n=d+1}. 
\end{proof}

%********************************************************************************

\section{The topological representation theorem --- general case}
\label{sec:FLgeneral}

This section is devoted to the inductive step in the proof of
Theorem~\ref{thm:FL}. 
Let $n> r>2$. Suppose that Theorem~\ref{thm:FL} holds for all 
hyperline sequences of rank $r$ with less than $n$ elements and for all
hyperline sequences of rank less than $r$. 
Thus, for any non-empty $R\subset E_n$ if the contraction $X/R$ (resp.\ the
deletion $X\setminus R$) is defined, then there is  an essentially
unique arrangement $\mathcal A(X/R)$ (resp.\ $\mathcal A(X\setminus R)$) of
oriented pseudospheres in $\p^{r-1-|R|}$ (resp.\ in
$\p^{r-1}$) of full rank with $X/R=X(\mathcal A(X/R))$ (resp.\ with
$X\setminus R = X(\mathcal A(X\setminus R))$).  

By Lemma~\ref{lem:fullrank}, there is an element of $X$, say, $n$ for
simplicity, so that $X\setminus \{n\}$ is a hyperline sequence of rank
$r$.  
Denote  $\{S_1,\dots,S_{n-1}\} = \mathcal A(X\setminus \{n\})$. Our aim 
is to construct an oriented pseudosphere $S_n\subset \p^{r-1}$ as the
image of a tame embedding $\psi\co \p^{r-2}\to \p^{r-1}$, so that 
$\{S_1,\dots,S_n\}$ is an arrangement of oriented pseudospheres with
$X(\{S_1,\dots, S_n\})= X$.

We outline informally the idea of the construction of $\psi$. We start
with the arrangement $\mathcal A(X/\{n\})$ in $\p^{r-2}$. We require
that $\psi$ maps this arrangement ``consistently'' to the arrangement
$\mathcal A(X\setminus \{n\})$, in the sense that any cell in $\mathcal
C\big(I,\mathcal A(X/\{n\})\big)$ is mapped to a cell in $\mathcal
C\big(I,\mathcal A(X\setminus\{n\})\big)$ in the correct orientation. It
turns out that this forces $\{S_1,\dots, S_n\}$ to be an arrangement 
of oriented pseudospheres.
Moreover, we show that if $S_n$ intersects the cycles of $\mathcal A(X
\setminus \{n\})$ in a way consistent with the rank 2 contractions of
$X$ (i.e., the cyclic order on its hyperlines),
then $X(\{S_1,\dots, S_n\}) = X$.
Our construction of $\psi$ is iterative. We start with defining $\psi$ on
$0$--dimensional cells of $\mathcal A(X/\{n\})$ and show that if it is
defined on $d$--dimensional cells then it can be consistently extended to
$(d+1)$--dimensional cells. It turns out that this is possible in an
essentially unique way.

Let us formalize this idea. 
By induction hypothesis, the arrangement $\mathcal A(X/\{n\})$ exists 
and is unique up to equivalence. For any element $i\in E(X/\{n\})$ of
$X/\{n\}$, we denote by $s_i$ the oriented pseudosphere of $\mathcal
A(X/\{n\})$ that corresponds to $i$.
For any $R\subset E(X/\{n\})$, set $s_R = \p^{r-2}\cap\bigcap_{j\in R}
s_j$, and similarly $S_R = \p^{r-1}\cap \bigcap_{j\in R} S_j$ for
$R\subset E_{n-1}$. 
Recall that the $d$--dimensional skeleton
$\mathcal A^{(d)}$ of an arrangement $\mathcal A$ of oriented
pseudospheres is the union of its cells of dimension $\le d$.
For any $R\subset E_n$, define $R/n = R\cap E(X/\{n\})$.

\begin{definition}
  Let $t< r$. A \begriff{$t$--admissible} embedding  is an embedding
  $$\psi^{(t)}\co \big(\mathcal A(X/\{n\})\big)^{(t-1)} \to \p^{r-1}$$
  so that  for any $R\subset E_{n-1}$ with
  $\dim s_{R/n} \le t-1$ holds 
  \begin{enumerate}
  \item $\psi^{(t)}(s_{R/n})=S_R$ or $\psi^{(t)}(s_{R/n})$ is a
    pseudosphere in $S_R$, 
  \item if $\psi^{(t)}(s_{R/n})\not=S_R$, then any cycle of $\mathcal
    A(X \setminus \{n\})$ in $S_R$ is either contained in
    $\psi^{(t)}(s_{R/n})$ or meets both connected components of
    $S_R\setminus \psi(s_{R/n})$, and
  \item for any $i\in E(X/\{n\})\setminus R\;$ holds 
    $\psi(s_{R/n}\cap s_i^+)  \subset S_R\cap S_i^+$ and
    $\psi(s_{R/n}\cap s_i^-) \subset S_R\cap S_i^-$.
  \end{enumerate}
\end{definition}

By the following two lemmas, in our request for the pseudosphere $S_n$ it
suffices to study $(r-1)$--admissible embeddings.
\begin{lemma}\label{lem:admissibelarrangement}
  Let $\psi\co\p^{r-2}\to\p^{r-1}$ be a tame embedding that defines an
  oriented pseudosphere $S_n$. If $\psi$ is $(r-1)$--admissible
  then $\mathcal A= \{S_1,\dots,S_n\}$ is an arrangement of oriented
  pseudospheres.
\end{lemma}
\begin{proof}
  We prove that $\mathcal A = \{S_1,\dots,S_n\}$ satisfies
  Axioms~(A1) and~(A2).
  \begin{itemize}
  \item[(A1)] Let $R\subset E_n$. It is to show that $S_R$ is empty or
    homeomorphic  to a sphere. If $n\not\in R$ then we are done since
    $\{S_1,\dots, S_{n-1}\}$ is an arrangement. 
    If $n\in R$ then we have $S_R=\psi(s_{R/n})$, and $s_{R/n}$ is empty 
    or homeomorphic to a
    sphere since $\mathcal A(X/\{n\})$ is an arrangement. 
  \item[(A2)] Let $R\subset E_n$ and $i\in E_n$ with $S_R\not\supset
    S_i$. Since $\psi$ is $(r-1)$--admissible, $S_R\cap S_i$ is a
    pseudosphere in $S_R$.  It remains to show that $S_R\cap S_i^+$ and
    $S_R\cap S_i^-$ are both non-empty.  If $R\cup \{i\}\subset E_{n-1}$
    then we are done, since $\{S_1,\dots, S_{n-1}\}$ is an arrangement
    of oriented pseudospheres.

    If $n\in R$ then $S_R= \psi(s_{R/n})$ and $S_R\cap S_i =
    \psi(s_{(R\cup \{i\})/n})$. Since $S_R\not\supset S_i$, we have
    $s_{R/n}\not\supset s_i$. Thus in this case, Axiom~(A2) for
    $\mathcal A$ follows from Axiom~(A2) for $\mathcal A(X/\{n\})$
    applied to $s_{R/n}$ and $s_i$, since
    $\psi(s_{R/n}\cap s_i^+) \subset S_R\cap S_i^+$ and  
    $\psi(s_{R/n}\cap s_i^-) \subset S_R\cap S_i^-$.

    If $n=i$, then $R=R/n$ since otherwise $S_n=S_i\subset S_R$ by
    definition of $S_n$.  Hence  $S_R\cap S_n = \psi(s_{R})$ is a
    pseudosphere in $S_R$. Since $\mathcal A\setminus \{n\}$ is of full
    rank, $S_R$ contains a cycle of $\mathcal A\setminus \{n\}$. By the
    second property in the definition of $(r-1)$--admissible embeddings,
    applied to the empty set, 
    this circle meets both  connected components of $\p^{r-1}\setminus
    S_n$, which implies Axiom~(A2). 
  \end{itemize}
  Thus, $\mathcal A$ is an arrangement of oriented
  pseudospheres. 
\end{proof}

\begin{lemma}\label{lem:arrangementadmissibel}
  Let $\tilde\psi\co\p^{r-2}\to\p^{r-1}$ be a tame embedding that defines an
  oriented pseudosphere $\tilde S_n$. Assume that $\tilde{\mathcal
  A}=\{S_1,\dots,S_{n-1},\tilde S_n\}$ is an arrangement of oriented 
  pseudospheres.
  If $X(\tilde{\mathcal A}) = X$ then there is an orientation preserving
  homeomorphism $\phi\co \p^{r-2}\to 
  \p^{r-2}$ so that $\tilde\psi\circ \phi$ is $(r-1)$-admissible.
\end{lemma}
\begin{proof}
  The arrangement $\tilde{\mathcal
    A}/\{n\}$ of oriented pseudospheres is given by the pre-images of
  $S_1,\dots,S_{n-1}$ under $\tilde\psi$.  By induction hypothesis in
  the proof of Theorem~\ref{thm:FL}, the equivalence class of ${\mathcal
    A}(X/\{n\})$ is uniquely determined by the property $X({\mathcal
    A}(X/\{n\})) = X/\{n\}$.  If $X(\tilde{\mathcal A}) = X$, then
  $X(\tilde{\mathcal A}/\{n\}) = X/\{n\}$. Hence, the arrangement 
  ${\mathcal A}(X/\{n\})$ is equivalent to $\tilde{\mathcal A}/\{n\}$.
  Let $\phi\co \p^{r-2}\to \p^{r-2}$ be the orientation preserving
  homeomorphism realizing this equivalence.  It easily follows  that
  $\psi\circ\phi$ is $(r-1)$--admissible.
\end{proof}

According to the preceding two lemmas, we shall construct $S_n$ via
a tame $(r-1)$--admissible embedding. Moreover, we need to take into
account the rank $2$ contractions of $X$, as follows.
Let $\psi^{(1)}$ be an $1$--admissible embedding. 
Assume that for any contraction $X/R$ of rank $2$ with $n\in E(X/R)$, the
oriented $0$--dimensional sphere $\psi(s_{R})$ extends the arrangement of
oriented pseudospheres on the oriented cycle $S_R$ that is induced by
$\mathcal A(X\setminus\{n\})$ to an arrangement equivalent to $\mathcal
A(X/R)$.   
Then we call $\psi^{(1)}$ \begriff{compatible with $X$}. 
In the next three lemmas, we prove that there is an essentially unique
$(r-1)$--admissible embedding whose restriction to the $0$--skeleton of
$\mathcal A(X/\{n\})$ is compatible with $X$. 

\begin{lemma}\label{lem:1admissibel}
  There is a tame $1$--admissible embedding $\psi^{(1)}$ that is
  compatible with $X$. It is unique up to composition with a
  homeomorphism $\big(\mathcal A(X\setminus \{n\})\big)^{(1)}\to \big(\mathcal
  A(X\setminus \{n\})\big)^{(1)}$ that fixes $\big(\mathcal
  A(X\setminus \{n\})\big)^{(0)}$. 
\end{lemma}
\begin{proof}
  Let $R\subset E_{n-1}$ so that $X/R$ is of rank $2$.
  We first prove the uniqueness of $\psi^{(1)}$.

  If $s_{R/n}\approx \p^0$, then the cyclic order of the signed elements of
  $X/R$ uniquely determines in which
  cells of $\mathcal A(X\setminus \{n\})$ on $S_R$ the two points of $s_R$ must
  be mapped to, provided $\psi^{(1)}$ is compatible with $X$. If 
  they are mapped to $0$--dimensional cells of $\mathcal
  A(X\setminus\{n\})$ then their image is unique. If they are mapped
  to $1$--dimensional cells, then their image is unique up to a
  homeomorphism of these cells fixing the boundary.

  If $\dim s_{R/n}>0$, then $R\not= R/n$, and $s_{R/n}$ is a  cycle of
  $\mathcal A(X/\{n\})$.  
  Any $0$--cell on $s_{R/n}$ is contained in some pseudosphere $s_i$ of
  $\mathcal A(X/\{n\})$ that intersects $s_{R/n}$ transversely. 
  Then $s_{R/n}\cap s_i\approx \p^0$, and
  if $\psi^{(1)}$ is $1$--admissible then $\psi(s_{R/n}\cap s_i)\subset
  S_R\cap S_i\approx \p^0$. Moreover, if $\psi^{(1)}$ is $1$--admissible
  then the intersection of $s_{R/n}\cap s_i$ with one side of a
  pseudosphere $s_j$ is mapped to the corresponding side of $S_j$. This
  uniquely determines the image of the two points of $s_{R/n}\cap s_i$ under
  $\psi^{(1)}$.

  We prove the existence of $\psi^{(1)}$. According to the preceding two
  paragraphs, for any $0$--dimensional cell $p$ of $\mathcal A(X/\{n\})$,
  a candidate for $\psi^{(1)}(p)$ is given by the cyclic order of the rank
  $2$ contractions of $X$.  
  We must ensure that the candidate does not depend on the choice of the
  contraction.
  Since any two cycles of $\mathcal A(X\setminus \{n\})$ are contained
  in some $2$--sphere that is the intersection of pseudospheres in
  $\mathcal A(X\setminus \{n\})$, we can assume by a contraction that
  $X$ is of rank $3$. Then, the pseudospheres $S_1,\dots,S_{n-1}$ are
  cycles. 
  
  Let $i,j\in E_{n-1}$. If $i\not\in E(X/\{n\})$ then $X/\{i\} = \pm
  X/\{n\}$, hence, the cyclic order of signed points on $S_i$ is a copy
  of $\mathcal A(X/\{n\})$, possibly with opposite orientation. In this
  case (and similarly if $j\not\in E(X/\{n\})$) it is easy to show that
  the candidates for $\psi^{(1)}$ imposed by $i$ and $j$ coincide. 
 It remains the case $i,j\in E(X/\{n\})$ with $s_i\cap s_j\not=\emptyset$.
 The positive point $p$ of the oriented pseudosphere $s_i$ shall be mapped
 into the cell $C_i$ of  $\mathcal A(X\setminus \{n\})$ on $S_i$ that
 corresponds  to the atom of $X/\{i\}$ containing $\overline n$. 
 If $(X/\{n\})\big/\{i\} = (X/\{n\})\big/\{j\}$
 (resp. $(X/\{n\})\big/\{i\} = -(X/\{n\})\big/\{j\}$), then $p$ is the
 positive (resp.\ negative) point of $s_j$, thus shall be
  mapped into the cell $C_j$ on $S_j$ that corresponds to the atom of
  $X/\{j\}$ containing $\overline n$ (resp.\  containing $n$). 

 If $X/\{i\} = \pm X/\{j\}$ then obviously $C_i=C_j$. Otherwise,
 $S_i\cap S_j$ is a $0$--sphere containing both $C_i$  and $C_j$. Up to
 symmetry, we can assume that $(X/\{i\})\big/\{n\}=(X/\{i\})\big/\{j\}$,
 which means that the atom of $X/\{i\}$ containing $\overline n$ does
 also contain $\overline \jmath$. Hence, $C_i$ corresponds to the atom of $X/\{j\}$
 containing $i$. This atom also contains $\overline n$ (resp.\ $n$) if
 and only if  
 $$  (X/\{j\})\big/\{n\} = - (X/\{j\})\big/\{i\} = (X/\{i\})\big/\{j\} =
 X/\{i\})\big/\{n\}$$
 (resp. (X/\{j\})\big/\{n\} = - (X/\{i\})\big/\{n\}). Thus $C_i= C_j$ by
 construction of $C_i$ and $C_j$. 
 In conclusion, the candidates for $\psi^{(1)}(p)$ imposed by $i$ and $j$
 coincide, which is enough to prove the existence of $\psi^{(1)}$
\end{proof}

\begin{lemma}\label{lem:extendtraces}
  Let $t<r - 1$, and let $\psi^{(t)}$ be a $t$--admissible
  embedding. 
  For any $t$--dimensional cell $c$ of $\mathcal A(X/\{n\})$
  there is a cell $c'$ of $\mathcal A(X\setminus\{n\})$ of dimension
  $t$ or $t+1$ so that $\psi^{(t)}(\d c) \subset \d c'$.
\end{lemma}
\begin{proof}
  Let $R\subset E_{n-1}$ be maximal so that $c\subset s_{R/n}$. In
  particular, $\dim s_{R/n}=\dim c= t$ and $\dim S_R\le \dim s_{R/n} +1=
  t+1$. For any cell $b$ of dimension $\dim b=\dim c -1$ in $\d c$ and
  any $j_b\in E(X/\{n\})$ with $b\subset s_{R/n}\cap s_{j_b}$, we have
  $\psi^{(t)}(b)\subset S_R\cap S_{j_b}$ since $\psi^{(t)}$ is
  $t$--admissible. By consequence, $\psi^{(t)}(\d c)\subset S_R$.
  
  Since $c$ is a cell of $\mathcal A(X/\{n\})$, $\d c\cap
  s_j^+=\emptyset$ or $\d c\cap s_j^-=\emptyset$, for all $j\in
  E(X/\{n\})$.
  Since $\psi^{(t)}$ is $t$--admissible, $\psi^{(t)}(\d c)\cap
  S^+_j=\emptyset$ or $\psi^{(t)}(\d c)\cap S^-_j=\emptyset$, for all
  $j\in E_{n-1}$. Thus, $\psi^{(t)}(\d c)$ is contained in the closure
  of a connected component $c'$ of $$S_R\setminus \bigcup_{j\in
    E_{n-1}\setminus R} S_j$$ which is a cell of $\mathcal
  A(X\setminus\{n\})$ of dimension $\dim S_R\le t+1$.
\end{proof}

\begin{theorem}\label{thm:psi_exists}
  There is an $(r-1)$--admissible embedding $\psi$ whose restriction to
  $\big(\mathcal A(X/\{n\})\big)^0$ is compatible with $X$. It is
  unique, up to a homeomorphism $\p^{r-1}\to \p^{r-1}$ that fixes
  $\mathcal A(X\setminus\{n\})$ cellwise.
\end{theorem}
\begin{proof}
  Let $\psi^{(t)}$ denote the restriction of $\psi$ to $\big(\mathcal
  A(X/\{n\})\big)^{(t-1)}$, for $t=1,\dots, r-1$.  
  We start with inductively proving the uniqueness of $\psi$.
  An $1$--admissible embedding $\psi^{(1)}$ that is compatible with $X$
  is essentially unique, by Lemma~\ref{lem:1admissibel}.
  For $t<r-1$, assume that $\psi^{(t)}$ is unique, up to a
  homeomorphism $\p^{r-1}\to \p^{r-1}$ that
  fixes $\mathcal A(X\setminus \{n\})$  cellwise.
  Let $c$ be a $t$--dimensional cell of $\mathcal A(X/\{n\})$, and let
  $R\subset E_{n-1}$ be maximal so that $\d c\subset s_{R/n}$.
  If $\psi$ is $(r-1)$--admissible, then $\psi(c)\subset S_R$. Let 
  $c'\subset S_R$ be the cell of $\mathcal A(X\setminus \{n\})$ containing
  $\psi(c)$. Since $\mathcal A(X\setminus \{n\})$
  is of full rank, the cell $c'$ is uniquely determined by $\psi(\d c)$.
  If $\dim c'=\dim c$ then the $(r-1)$--admissibility of $\psi$ imposes
  $c'=\psi (c)$.  Otherwise, it is a consequence of the generalized
  Sch{\"o}nflies theorem~\cite{brown} that $\psi(c)\subset c'$, being a
  tame ball of codimension $1$, is unique up to a homeomorphism
  $\overline c'\to \overline c'$ that fixes $\d c'$ pointwise.
  This can be extended to a homeomorphism $\p^{r-1}\to \p^{r-1}$ fixing
  all cells of $\mathcal A(X\setminus \{n\})$, by the so-called cone
  construction~\cite{stoeckerzieschang}.
  In conclusion, $\psi^{(t+1)}$ is essentially uniquely determined by $\psi^{(t)}$,
  which completes the proof of the uniqueness of $\psi$.
  
  Secondly, we expose an iterative construction of $\psi$.  An
  $1$--admissible embedding $\psi^{(1)}$ that is compatible with $X$
  exists, by Lemma~\ref{lem:1admissibel}.
  For $t<r-1$, assume that $\psi^{(t)}$ is a $t$--admissible embedding,
  and let $c$ be a $t$--dimensional cell of $\mathcal A(X/\{n\})$.  If
  there is a $t$--dimensional cell $c'$ of $\mathcal A(X\setminus
  \{n\})$ so that $\d c'=\psi^{(t)}(\d c)$ then we define
  $\psi^{(t+1)}(c) = c'$.
  Otherwise, by the preceding lemma there is a $(t+1)$--dimensional cell
  $c'$ of $\mathcal A(X\setminus \{n\})$ so that $\psi^{(t)}(\d
  c)\subset \d c'$. 
  Since $\psi^{(t)}(\d c)$ is a tame union of tame cells of codimension
  one in $\d c'$, it follows from~\cite{kirby} that $\psi^{(t)}(\d c)$
  is a pseudosphere in $\d c'$.  Therefore $\psi^{(t)}(\d c)$ bounds a
  tame cell $c''\subset c'$.  We define $\psi^{(t+1)}(c) = c''$.
  
  We show that we can do this construction so that $\psi^{(t+1)}$ is an
  embedding.  Let $\Psi^{(t)}$ denote the image of $\psi^{(t)}$, and let
  $c'$ be a $(t+1)$--dimensional cell of $\mathcal A(X\setminus \{n\})$.
  If $\psi^{(t)}$ is an embedding, then $\Psi^{(t)}\cap \d c'$ is a disjoint
  union of $t$--dimensional spheres. If $t>1$ then these spheres bound
  a system of disjoint $(t+1)$--dimensional cells in $c'$. If $t=1$ then
  a priori the spheres might be linked, as depicted in
  Figure~\ref{fig:linked}. The thick dots indicate the images under $\psi$ of
  four $0$--dimensional cells on a cycle of $\mathcal A(X/\{n\})$, 
  whose cyclic order corresponds to the numbering.
  \begin{figure}[htbp]
    \begin{center}
      \epsfig{file=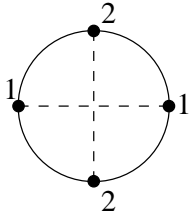}
      \caption{Two linked $0$--spheres in $\p^1$}
      \label{fig:linked}
    \end{center}
  \end{figure}
  It is shown in~\cite{bokowskimockstreinu} 
  that this case does not occur.
  In conclusion, our construction of $\psi^{(t+1)}$ produces an
  embedding.

  It remains to show that $\psi^{(t+1)}$ is $(t+1)$--admissible. Let
  $R\subset E_{n-1}$ with $\dim s_{R/n} \le t$.
  \begin{enumerate}
  \item Either $\psi^{(t+1)}(s_{R/n})$ is equal to $S_R$ or it is
    composed by tame cells of codimension one in $S_R$. In the latter
    case, $\psi^{(t+1)}(s_{R/n})$ is a pseudosphere (i.e., tame)
    by~\cite{kirby}.
  \item Let $\psi^{(t+1)}(s_{R/n})\not=S_R$, and let $\gamma\subset S_R$
    be a cycle of $\mathcal A(X\setminus \{n\})$ that is not contained
    in $\psi^{(t+1)}(s_{R/n})$. 
    Chose a subset $T\subset E_{n-1}$ with $\gamma = S_T$ and $|T|=
    r-2$.  
    The contraction $X/T$ is a hyperline sequence of rank $2$. 
    Since $\psi^{(1)}$ is compatible with $X$,
    $\psi^{(t+1)}(s_{R/n})\cap \gamma$ comprizes exactly two points
    $x,y$, corresponding to the elements $n,\overline n$ in $X/T$.
    
    There is some $x\in\mathbf E(X/T)\setminus \{n\}$, so that $X/T$
    induces the cyclic order $(n,x,\overline n, \overline x)$. Let
    $j=x^*$. Since $\psi^{(t)}$ is $t$--admissible, it follows from
    Lemma~\ref{lem:admissibelarrangement} that $\mathcal A(X\setminus
    \{n\})$ and $\psi^{(t+1)}(s_{R/n})$ induce on $S_{R\cup\{j\}}$ an
    arrangement of oriented pseudospheres. If $\gamma$ does not
    meet both connected components of $S_R\setminus
    \psi^{(t+1)}(s_{R/n})$, then the two points of $\gamma\cap S_j$ are
    contained in a single component of $S_{R\cup\{j\}}\setminus
    \psi^{(t+1)}(s_{R/n})$, which is impossible for arrangements of
    oriented pseudospheres. Hence, $\psi^{(t+1)}$ satisfies the second
    property in the definition of $(t+1)$--admissible embeddings.
  \item We observe that an open cell is contained in the positive side of an
    oriented pseu\-do\-sphe\-re if and only if some point in its boundary is
    contained in the positive side of the oriented pseudosphere.  Thus,
    if $c$ is a $t$--dimensional cell of $\mathcal A(X/\{n\})$ and
    $c\subset s_i^+$, then $\d c \cap s_i^+\not= \emptyset$. Since
    $\psi^{(t)}$ maps $s_i^+$ into $S_i^+$, it follows
    $$\psi^{(t)}(\d c) \cap S_i^+ = \psi^{(t+1)} (\d c)\cap S_i^+ \not=
    \emptyset,$$
    and therefore $\psi^{(t+1)}(c)\subset S_i^+$.
    Similarly, if $c\subset S_i^-$ then $\psi^{(t+1)}(c)\subset S_i^-$.
  \end{enumerate}
  Therefore $\psi^{(t+1)}$ is $(t+1)$--admissible, which finishes the
  proof of Theorem~\ref{thm:psi_exists}.
\end{proof}

In the remainder of this section, we finish the proof of
Theorem~\ref{thm:FL}. 
Let $\psi\co \p^{r-2}\to \p^{r-1}$ be a $(r-1)$--admissible embedding
whose restriction to the $1$--skeleton of $\mathcal A(X/\{n\})$ is
compatible with $X$. By Theorem~\ref{thm:psi_exists}, $\psi$ exists. 
Let $S_n=\psi(\p^{r-2})$ be the oriented  pseudosphere defined by $\psi$,
and set $\mathcal A = \{S_1,\dots,S_n\}$. 
By Lemma~\ref{lem:admissibelarrangement}, $\mathcal A$ is an arrangement
of oriented pseudospheres. 

We show that $X(\mathcal A) = X$, hence, that the hyperlines $(Y|Z)\in
X$ coincide with those of $X(\mathcal A)$.
If $n\not\in E(Y)$ then $(Y|Z)\in X(\mathcal A)$, since $Z$ is a rank
$2$ contraction of $X$ and $\psi^{(1)}$ is $1$-admissible and compatible
with $X$.  
If $n\in E(Y)$, then let $R=E(Y)\setminus \{n\}$. Since $\psi$ is
$(r-1)$--admissible, we have $S_R\subset S_n$, thus, $S_R$ corresponds
to the cycle $s_{R/n}\subset \p^{r-2}$ of $\mathcal A(X/\{n\})$. 
Since $X\big(\mathcal A(X/\{n\})\big)= X/\{n\}$, the cyclic order of
points on this cycle coincides with $Z$. If $i,j\in E(Z)$ so that $Z$
induces the cyclic order $(i,j,\overline \imath,\overline \jmath)$, then
$X/\{i,j\}=Y$. We have $X(\mathcal A/\{i,j\}) = X/\{i,j\}$, since 
the Topological Representation Theorem~\ref{thm:FL} holds in rank
$r-2$ by induction hypothesis. Therefore, $(Y|Z) = (X(\mathcal
A/\{i,j\})|Z)\in X(\mathcal A)$. 

It remains to prove the uniqueness of $\mathcal A$ stated in the
Topological Representation Theorem~\ref{thm:FL}.  Let $\{\tilde
S_1,\dots, \tilde S_n\}$ be an arrangement of oriented pseudospheres
with $X(\{\tilde S_1,\dots, \tilde S_n\})=X$. By induction hypothesis,
we may use the uniqueness of $\mathcal A(X\setminus \{n\})$ and can
assume that $\{\tilde S_1,\dots,\tilde S_{n-1}\} = \{S_1,\dots,
S_{n-1}\}$. Let $\tilde S_n$ be the image of a tame embedding
$\tilde\psi\co\p^{r-2}\to \p^{r-1}$.  Then, $\tilde\psi$ is
$(r-1)$--admissible by Lemma~\ref{lem:arrangementadmissibel}, and it is
obvious that its restriction to the $1$-skeleton of $\mathcal
A(X/\{n\})$ is compatible with $X$.  Thus by Theorem~\ref{thm:psi_exists},
$\tilde \psi$ coincides with $\psi$ up to a homeomorphism $\p^{r-1}\to
\p^{r-1}$ that fixes $\{S_1,\dots,S_{n-1}\}$ cellwise.  Hence, $\mathcal A$
is equivalent to $\{\tilde S_1,\dots, \tilde S_n\}$.  This finishes the
proof of the Topological Representation Theorem~\ref{thm:FL}.

%***************************************************************************

\section{Wild arrangements}
\label{sec:wild}

This section is an appendix to Section~\ref{sec:X(A)}. We show here that
one can get a hyperline sequence from an arrangement of oriented
pseudospheres, even if one allows pseudospheres to be not tame.

Let $S_1,\dots, S_n\subset \p^d$ be embedded $(d-1)$--dimensional
spheres with a choice of an orientation.  We call the ordered multiset
$\mathcal A = \{S_1,\dots, S_n\}$ an \begriff{arrangement of oriented
  topological spheres} over $E_n$ if it satisfies Axioms~(A1) and~(A2),
where the word ``pseudosphere'' is replaced by ``embedded sphere of
codimension one''.

An embedded sphere $S\subset \p^d$ of codimension one is \begriff{wild},
if there is no homeomorphism $\p^d\to \p^d$ mapping $S$ to
$\p^{d-1}\subset \p^d$.  Since there are infinitely many wild spheres in
$\p^d$ for all $d\ge 3$ (see~\cite{rolfsen}), we are no longer allowed
to use the generalized Sch\"onflies theorem~\cite{brown}.  However,
$\p^d\setminus S$ has exactly two connected components, even if $S$ is
wild. 

Similarly to Section~\ref{sec:PSA}, we say that two arrangements
$\{S_1,\dots, S_n\}$ and $\{\tilde S_1,\dots, \tilde S_n\}$ of oriented
topological spheres in $\p^d$ are equivalent if there is an orientation
preserving homeomorphism $\phi\co \p^d \to \p^d$ with $\phi(S_i)=\tilde
S_i$ in the correct orientation, for $i=1,\dots, n$.
For $R\subset E_n$, the definition of the contraction $\mathcal A/R$ and
the deletion $\mathcal A\setminus R$ is identic to the corresponding
definition for arrangements of oriented \emph{pseudo}spheres.

We wish to define a hyperline sequence $X(\mathcal A)$ associated to an
arrangement $\mathcal A$ of oriented topological spheres.  By
Section~\ref{sec:X(A)}, we know how to proceed if all spheres in
$\mathcal A$ are tame. Both in the construction of $X(\mathcal A)$ in
Section~\ref{sec:X(A)} and in the proof that $X(\mathcal A)$ is indeed a
hyperline sequence, we were using induction on the contractions of
$\mathcal A$, based on contractions of rank $1$ and $2$. 
The only topological argument in the induction step was the use of the
Jordan--Sch\"onflies theorem, that holds also without the assumption of
tameness though, in the step from rank $2$ to rank $3$.
It remains to remove the tameness condition in the base cases.
Rank $1$ is trivial. 
For rank $2$, observe that any two different points $x,y\in \p^1$ can be
separated by small intervalls around $x$ and $y$. Thus, any embedded
sphere $\p^0$ in $\p^1$ is tame.
Therefore, even if $\mathcal A$ is not equivalent to an arrangement of
pseudospheres, any rank $2$ contraction $\mathcal A/R$ actually is an
arrangement of pseudospheres, and we can read off a hyperline sequence
$X(\mathcal A/R)$. As in Section~\ref{sec:X(A)}, these contractions
yield a hyperline sequence $X(\mathcal A)$.

In conclusion, although arrangements of oriented topological spheres are very
complicated from a topological point of view, their combinatorics is
simple enough to read off ordinary hyperline sequences. Nevertheless,
there are \lq\lq more\rq\rq\ arrangements of oriented topological spheres than hyperline
sequences, in the sense that there are non-equivalent arrangements
$\mathcal A_1$ and $\mathcal A_2$ of oriented topological spheres (for instance, a
tame arrangement and a wild arrangement) with $X(\mathcal A_1) =
X(\mathcal A_2)$. Hence, there is no analogue of the Topological
Representation Theorem~\ref{thm:FL} in the setting of arrangements of
oriented topological spheres.\\

\noindent {\bf Acknowledgment:}
We wish to thank Gerhard Burde, Klaus Johannson and Martin Scharlemann
for their advises for topological references.\\ 

%\newpage
%{\small

%}


\begin{thebibliography}{99}
\addcontentsline{toc}{section}{References}


\bibitem{bjorner}
 A. Bj\"{o}rner, M. Las Vergnas, B. Sturmfels, N. White and G. Ziegler,
\newblock {\em Oriented Matroids},
\newblock Cambridge University Press, (1993) second edition 1999.

\bibitem{bokowski}
 J. Bokowski,
\newblock {\em Oriented Matroids},
\newblock in \cite{gruber} (1993) 555-602.

\bibitem{bokowskimockstreinu}
J. Bokowski, S. Mock, I. Streinu,
\newblock {\em On the Folkman-Lawrence topological representation theorem
for oriented matroids in rank 3.}
\newblock {\em European J. of Combinatorics, 22,5 (2001) 601-615.}

\bibitem{brown}
\newblock Brown, M.:
\newblock A proof of the generalized Schoenflies theorem.
\newblock {\em Bull. Amer. Math. Soc.,} {66,} (1960) 74-76.

\bibitem{brown2}
\newblock Brown, M.:
\newblock Locally flat imbeddings of topological manifolds. 
\newblock {\em Ann. of Math. (2),} {77,} (1962) 331--341.

\bibitem{edmonds}
 J. Edmonds and A. Mandel,
\newblock{\em Topology of Oriented Matroids},
\newblock Abstract 758-05-9, {\em Notices Amer.Math.Soc.} 25, A-410, 1978;
\newblock also Ph.D. Thesis of A. Mandel, Univ. of Waterloo, 1982.

\bibitem{folkman}  
 J. Folkman and J. Lawrence, 
\newblock {\em Oriented Matroids},
\newblock J. Comb. Th., Ser. B, 25, (1978) 199-236.

\bibitem{gp}
 J.E. Goodman and R.~Pollack,
\newblock {\em Semispaces of configurations, cell complexes of arrangements},
\newblock J. of Combinatorial Theory, Series A,
37 (1984) 257-293.

\bibitem{goodman}
 J.E. Goodman,
\newblock {\em Pseudoline arrangements},
\newblock in \cite{gor}, (1997) 83-110.

\bibitem{gor}
 J.E. Goodman and J.~O'Rourke (eds.),
\newblock {\em Handbook of Combinatorial and Computational Geometry},
\newblock CRC Press, 1997.

\bibitem{gruber}
 P.M. Gruber and J.M. Wills (eds.),
\newblock {\em Handbook of Convex Geometry},
\newblock vol. A and B, North Holland, 1993.

\bibitem{gruenbaum}
 B. Gr\"unbaum,
\newblock {\em Arrangements and Spreads},
\newblock Regional Conf. Ser. Math., Amer. Math. Soc., number 10,
Providence, RI, 1972.

\bibitem{GN}
 L. Gutierrez Novoa,
\newblock{\em On n-ordered sets and order completeness},
\newblock  Pacific J. Math., 15, (1965) 1337-1345.

\bibitem{hochstattler}
 W.~Hochst\"attler,
\newblock {\em Seitenfl\"achenverb\"ande orientierter Matroide},
\newblock PhD Thesis, K\"oln, 1992.

\bibitem{kirby}
  R.C.~Kirby,
  \newblock{\em The union of flat $(n-1)$--balls is flat in $R^n$},
  \newblock Bull. Amer. Math. Soc. 74, (1968) 614--617. 

\bibitem{knuth}
 D.E. Knuth,
\newblock {\em Axioms and Hulls},
\newblock Springer Verlag, New York, 1992.

\bibitem{lv1}
 M.~Las Vergnas,
\newblock {\em Matro\"{\i}des orientables},
\newblock C.R. Acad. Sci. Paris, Ser. A, 280, (1975) 61-64.

\bibitem{lv2}
 M.~Las Vergnas,
\newblock {\em Bases in oriented matroids},
\newblock J. Combinatorial Theory, Ser. B, 25, (1978) 283-289.

\bibitem{lawrence}
 J.~Lawrence,
\newblock {\em Oriented matroids and multiply ordered sets}, 
\newblock Linear Algebra Appl., 48, (1982) 1-12.

\bibitem{rgz}
 J.~Richter-Gebert and G.~Ziegler,
 \newblock{\em Oriented Matroids},
\newblock in \cite{gor}, (1997) 111-132.

\bibitem{rolfsen}
  D~Rolfsen,
  \newblock {\em Knots and Links},
  \newblock Mathematics Lecture Series 7.
  \newblock Publish or Perish (1976).

\bibitem{stoeckerzieschang}
R.~St\"ocker, H.~Zieschang,
\newblock \textit{Algebraische Topologie. Eine Einf\"uhrung.}
\newblock Second edition. 
\newblock B. G. Teubner, Stuttgart, 1994.

\end{thebibliography}
\end{document}